# Design, modeling and control of a hybrid grid-connected Photovoltaic-Wind system for the region of Adrar, Algeria

Y. Kebbati [1,2] and L. Baghli [3]

Email: yassine.kebbati@univ-evry.fr

[1] Department of Energy Engineering, Institute of Water and Energy Sciences (Incl. Climate Change), PAN African University, Tlemcen, Algeria

[2] Department of Electrical and Electronics Engineering, Paris Saclay University, Evry, France

[3] GREEN-ENSEM, University of Lorraine, Nancy, France

## ABSTRACT

The use of fossil energy for electricity production is an evident source of pollution, global warming and climate change. Consequently, researchers have been working to shift towards sustainable and clean energy by exploiting renewable an environmentally friendly resources such as wind and solar energies. On the other hand, energy security can only be achieved by considering multiple resources. Large-scale renewable energy power plants are a key solution for diversifying the total energy mix and ensuring energy security. This paper presents a contribution to diversify the energy mix in Algeria and help mitigate power shortages and improve grid performance. In particular, the paper aims at designing and modeling a large-scale hybrid Photovoltaic-wind system that is grid connected. An innovative control approach using improved particle swarm optimized PI controllers is proposed to control the hybrid system and generate the maximum power from the available wind and solar energy resources. Furthermore, economic, environmental and feasibility studies of the project were conducted using HOMER software to assess the viability of the system and its contribution to reduce greenhouse gas emissions.

**Keywords**: Renewable Energy Systems, Photovoltaic, Particle Swarm Optimization, Wind Turbine, Adaptive Control.

# 1. Introduction

In today's world with an ever-increasing energy demand due to the developing industry and modern lifestyle, fossil fuel resources are being overexploited. Such high consumption rates caused a significant surge in greenhouse gas (GHG) emissions which effectively lead to a rise in earth's temperature. The side effects of global warming are nowadays expressed more than ever leading to droughts, extreme weather, flooding, etc. In addition, reserves of oil and natural gas are being rapidly depleted. Dependance on fossil fuel has become a serious issue that needs urgent attention on any levels. For instance, the economic sector in Algeria relies mainly on hydrocarbons exportation. The country's daily production reaches 92.1 billion cubic meters of natural gas in addition to 1.54 million barrels of crude oil. This huge reliance on fossil fuels has a major effect on the economy of the country as the prices of oil and gas are rarely stable. Oil and gas exports account for nearly 20% of the gross domestic product, and about 85% of total exports (Mustapha et al. 2018). Over the last years, the fall of oil prices has widely affected the economy resulting in a huge deficit in the country's budget. Between 2014 and 2016, the price of one barrel of oil fell from \$100 to \$45.13 leading to a reduction of the total revenues by two thirds, dropping from \$60 billion in 2014 to less than \$20 billion in 2016. With this sharp decline, the country's gross domestic product (GDP) fell from \$235 billion in 2014 to \$160 billion in 2016. Moreover, fossil fuel prices dropped recently due to the COVID-19 pandemic and the global economic recession (Yoshino et al. 2020). In 2015 for instance, the state's budget reached its lowest record at -15.3% of the total GDP.

Electricity and oil prices are majorly subsidized in Algeria. For instance, electricity costs 0.04 \$/kWh which is among the lowest prices in the world (GlobalPetrolPrices.com 2020). However, since the emergence of the oil crisis, the subsidies started decreasing, and fossil fuel prices are now growing more and more expensive. Subsidized electricity prices lead to one of the most irrational consumptions in the continent. The latter reached 1302 kWh per capita in 2020 and it is expected to exceed 150 TWh in 2030 (Index Mundi 2021). In light of the projection studies and as fossil fuel resources are none renewable and unsustainable, Algeria's oil reserves are expected to only last for the next 50 years. Similarly, natural gas is expected to cover the next 70 years (Harrouz et al. 2017). The major dependance on fossil fuels and the irrational power consumption in Algeria results in considerable amounts of GHG emissions. Over the last decades, significant levels of GHG emissions have been reached, which causes an impactful deterioration of the ecological system increasing global warming and further accelerating climate change. There is no doubt that the shift towards renewable energy resources is a more ecological and sustainable choice. Thus, integrating renewable energies will reduce fossil energy dependency and dramatically decrease GHG emissions. Consequently, researchers have been working to further develop renewable energy systems in order to increase their efficiency and extract the maximum possible energy from renewable resources. The goal is to become independent of fossil fuel resources and cover 100% of power demand from renewable energies. For instance, (Amit et al. 2019) reviewed the maximum power point tracking techniques according to their 8 categories to increase the efficiency of photovoltaic systems. (Wijeratne et al. 2019) investigated the current approaches for the design and development of distributed solar PV systems and their efficacy. (Mahmoud et al. 2020) analyzed the effect of cold weather conditions on the efficiency of PV systems over the years, their studies were based on the weather patterns in UK and in Australia. (Gopal et al. 2021) provided a solution to increase the performance of wind turbine systems while reducing their overall cost. Their approach relies on minimizing the operating temperature of generator windings and bearings which reduces overheating of the system. (Kumar et al. 2021) carried out a technical review on the methods and techniques for damage detection in wind turbine systems using several signal processing methods. To enhance heat absorption from solar energy, (Djamel et al. 2021) conducted a study to use passive baffles in an attempt to improve the yield of single slope solar still. The latter is useful in water desalination applications. (Mustafa et al 2011) compared the economic viability of wind turbine and diesel water pumping systems in five different regions in Turkey, the authors explored various energy conversion systems to evaluate wind power and the amount of water to be pumped. (Gökçek et al. 2009) used the levelized cost of electricity to evaluate the effect of different wind energy conversion systems on the final wind energy cost. Similarly, (Fadlallah et al. 2021) used a strategic perspective to explore the viability of wind energy systems in Sudan. The authors used HOMER software to identify the optimal wind system and best locations for wind farm installation over 21 different sites. (Genç et al. 2021) provided a method to determine the most appropriate site for wind power systems. Their approach was based on geographic information combined with multi-criteria decision making that takes into account several factors like wind potential, roads and water sources. The same approach was used in (Genç et al. 2021) for site selection in the coastal line of turkey for offshore wind farms installations. Similarly (Shata et al. 2006) analyzed wind data for the Mediterranean see in Egypt to determine wind characteristics for assessing the viability of electricity generation using wind turbines.

Worldwide energy consumption in 2012 had a respectable share of 19% from renewable sources, the contribution of renewable energy continued to rise year after year reaching 19.3% in 2015, of which 1.6% goes to solar, wind, biomass and geothermal power combined (Corry et al. 2016). Algeria is home to one of the world's greatest solar deposits, the duration of sun exposure is more than 2,000 hours a year, up to 3900 hours in the highlands and the desert. This means that the country enjoys from 1,700 to 2,263 kWh/m$^2$/year of solar energy (Maoued et al. 2015). The south of Algeria has significant wind resources, especially the region of Adrar where average wind speeds range from 4 to 6 m/s, which makes it very attractive for the deployment of wind farms (Maoued et al. 2015). Fig. 1 shows the annual distribution of mean wind speed in Algeria at 10 meters altitude and the average daily solar irradiation (Boudia et al. 2016).

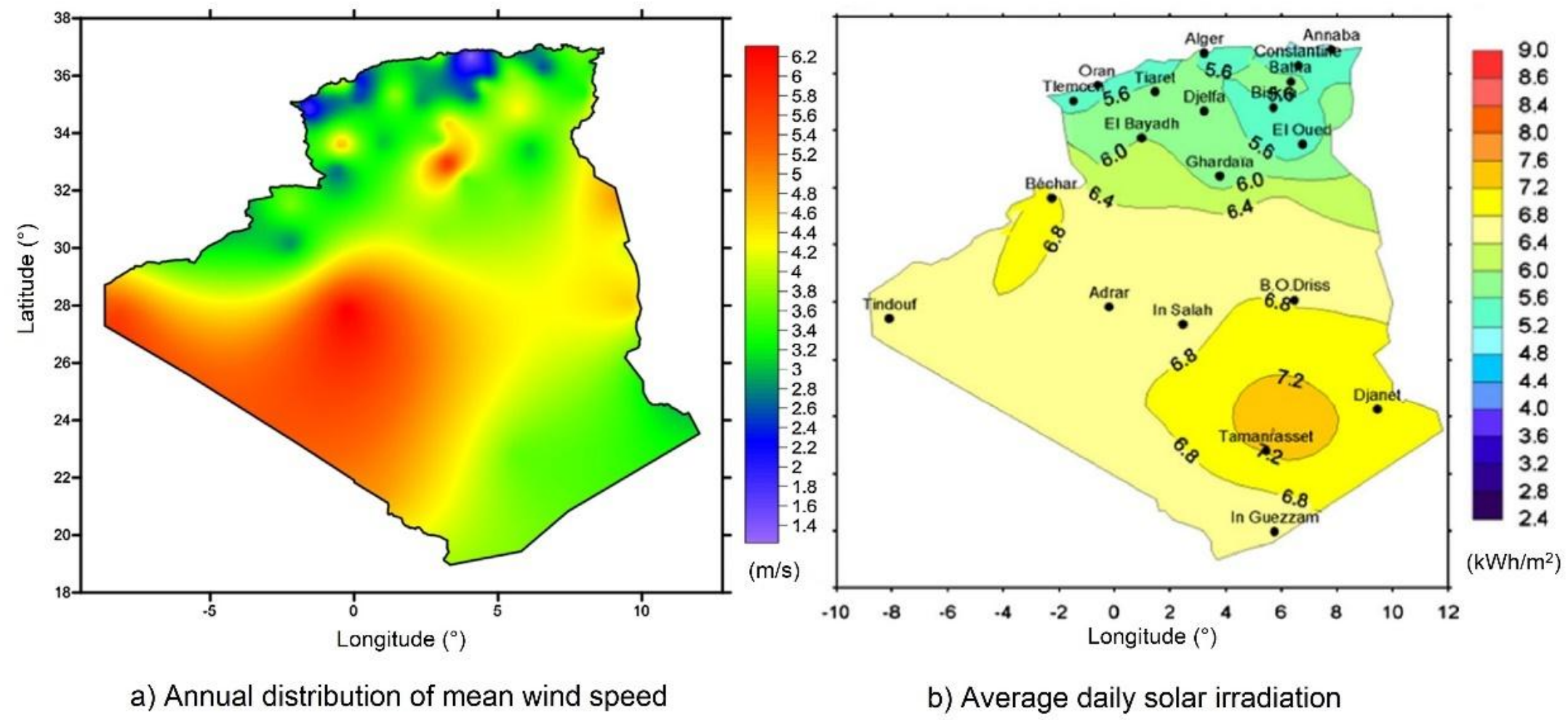

**Fig. 1: Mean wind speed and average solar irradiation in Algeria (Boudia et al. 2016).**

The Algerian government is trying to attract investments in wind and solar energies by establishing suitable policies to install 5 GW of wind power and 13.6 GW of solar PV by 2030. Moreover, the country launched in 2011 an ambitious renewable energy and energy efficiency program with the vision of developing and diversifying renewable energy resources. The goal of the program is to install 22 GW by 2030, of which 10 GW is aimed to be exported to the European market and 12 GW is dedicated to the national market (Abada et al. 2018). The distribution and share of renewable energy technologies is shown in Fig. 2.

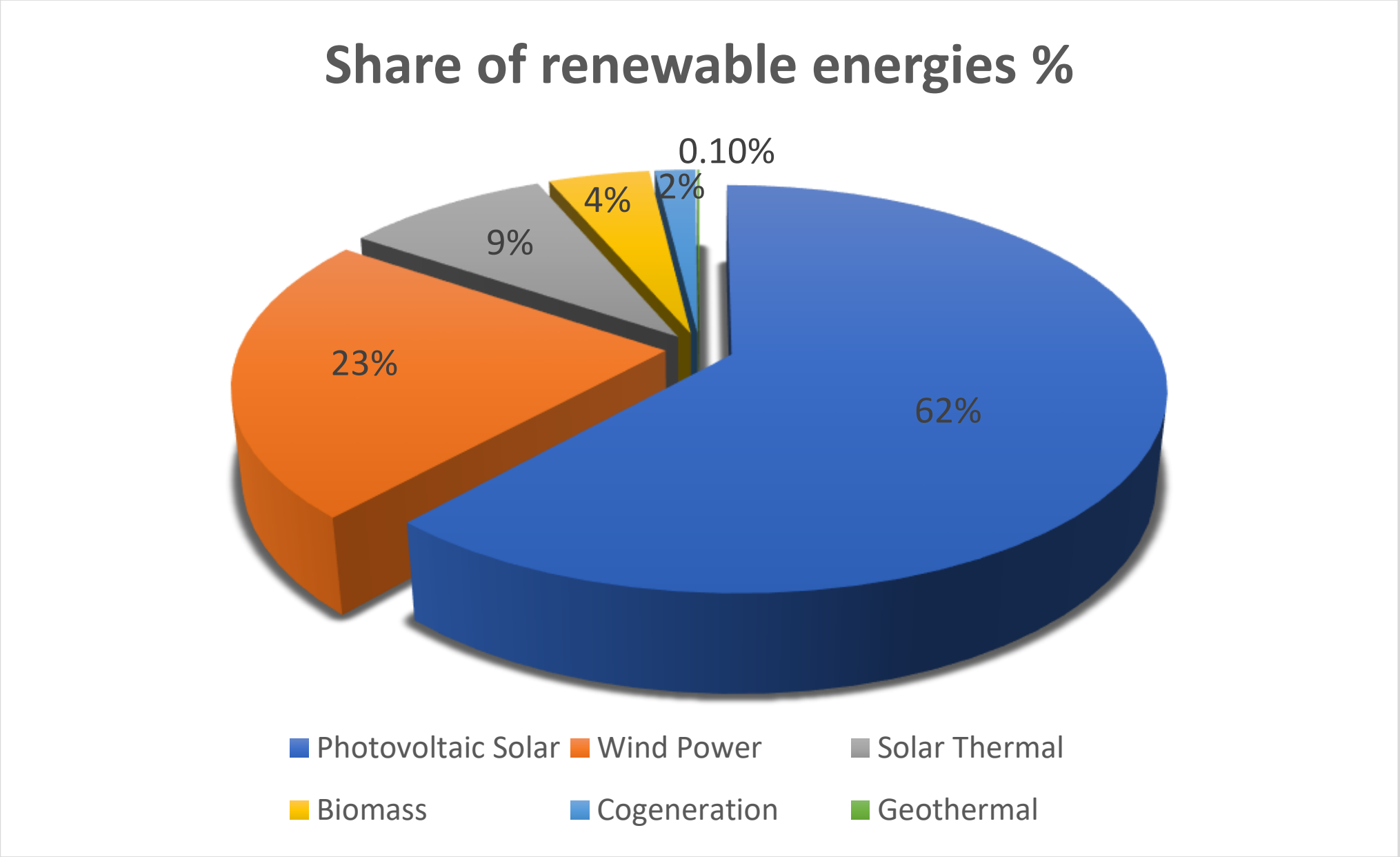

**Fig. 2: Share of renewable energy technologies in energy efficiency program.**

Algeria is known for its important potential in hydrocarbon resources, but it also houses a huge solar energy potential (Kabir et al. 2018; Bouraiou et al. 2020) thanks to its location in the Mediterranean basin. The southern part of the country receives more sunshine as it lies exactly in the sunbelt. Wind energy is the second most abundant resource; half of the country's surface shows a considerable average wind speed that favors wind energy exploitation (Nedjari et al. 2018; HIMRI et al. 2020). Judging by the available resources, the integration of solar photovoltaic (PV) and wind energy systems is a key solution to further diversify the energy mix in the country and increase energy security. Grid-connected PV systems have proven their potential and capacity to provide adequate energy to the grid, these systems have been studied in (Kidar et al. 2021) and (Zerglaine et al. 2021). In highly distributed populations, extending power grids may become too expensive. For instance, (Terfa et al. 2021) proposed the use of distributed micro-power plants as a solution to the expensive grid extension. On the other hand, the southwest region is very suitable for wind energy systems which have been studied in (WANG et al. 2018) and (Charrouf et al. 2018). Furthermore, (Dekali et al. 2021) worked on the design, modeling, and experimental build of a 1.5 kW low-cost wind energy conversion system based on current-controlled DC motors and tip speed ratio (TSR) based on maximum power point tracking (MPPT).

Adrar represents the second largest province in Algeria, situated in the heart of the desert with an area of 427368 km². Due to its considerable wind and solar energy potential, it has been selected as the case study. The district of Adrar is precisely situated at 27.52° north of latitude and 0.17° west of longitude occupying 663 km² of the whole province. Around 65000 inhabitants live mainly in the city of Adrar, even though the province is vast. The climate in Adrar is a Saharan arid one. As for the topography, Adrar is characterized by a relatively flat terrain with the highest point reaching 421 meters. The location of the case study region is given in Fig. 3. Renewable energy resources vary in nature, solar and wind proved to be intermittent and this fact makes them a bit tricky to harness; because their intermittency creates challenges to maintaining a continuous and reliable power supply. In grid-connected systems, it is even more difficult especially in the case of weak grids that are not able to handle the fluctuation of power generation when the amount of integration of solar or wind is important. Hybrid systems can tackle this issue, combining solar PV with wind is an attractive solution that provides reliable and economical renewable power generation.

In this article, a hybrid grid-connected PV-wind system is designed, modeled and controlled with optimized PI controllers. A new improved particle swarm optimization (PSO) algorithm was developed to optimize the gains of the PI controller. The proposed system was assessed for two critical months of the year to evaluate its contribution to grid performance. Hourly weather data averaged over the last 10 years have been used to simulate the system. Finally, the economic and environmental analysis of the system was carried out using Homer software.

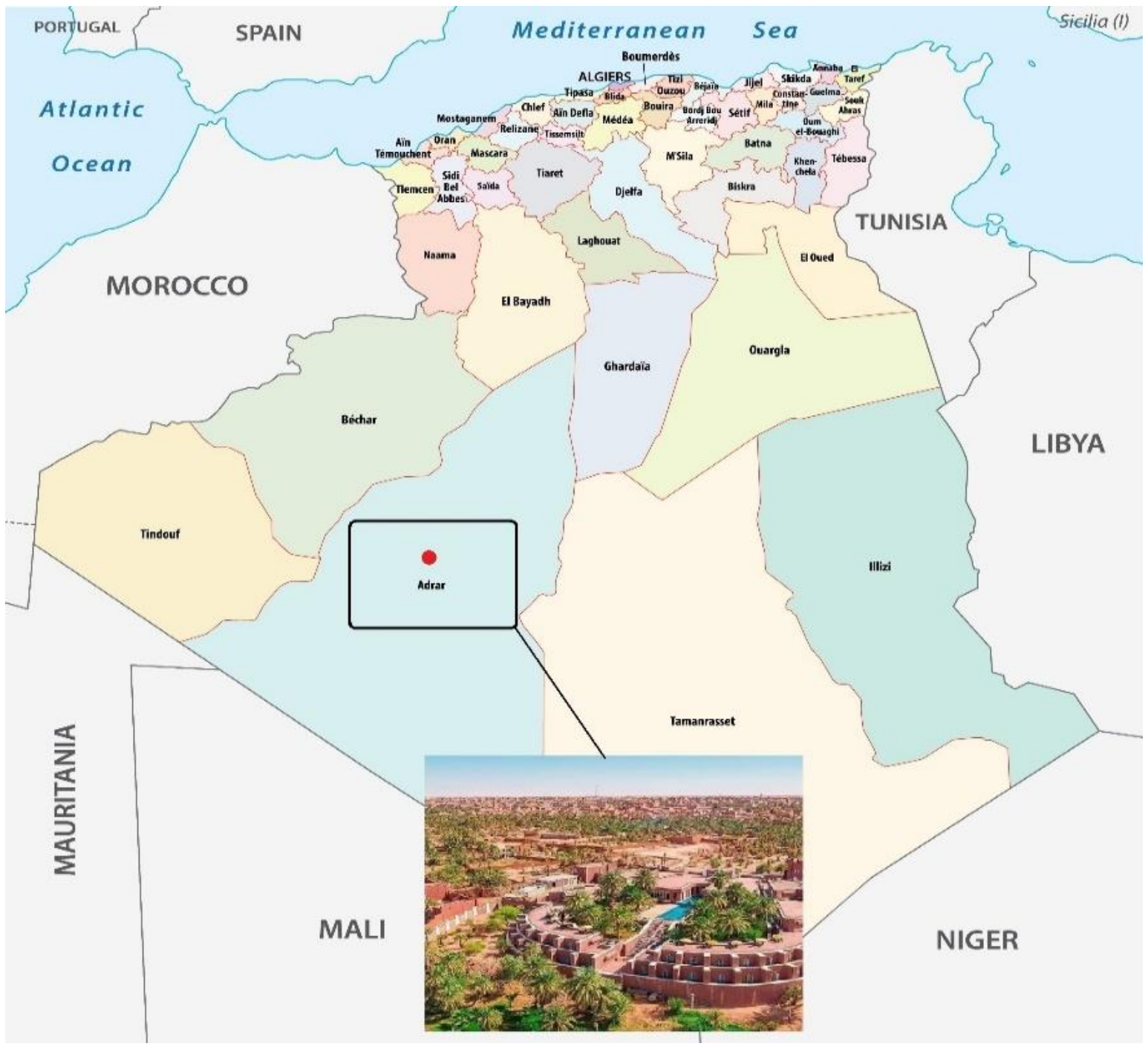

**Fig. 3: Location of Adrar on the map.**

# 2. Materials and Methods

## 2.1. Site selection and resource assessment

The province of Adrar has been selected for this study thanks to its potential in both wind and solar energies, it is situated in the heart of the Algerian desert located at 27.52° N and 0.17° W. The district occupies 663 km² and has a population of about 65000 inhabitants (Djamai et al. 2011). The climate in the province is a Saharan arid one with a topography characterized by flat terrain reaching the highest point of 421 meters. Based on hourly global irradiation data on a horizontal surface, the province of Adrar receives average daily irradiation of up to 5.7 kWh/m$^2$/day. During the winter season, the global radiation drops in between 3~4 kWh/m$^2$/day. The highest potential solar radiation is observed from March to October, with average global radiation varying between 5.5~7.5 kWh/m$^2$/day (Bouzidi B. 2011). The Algerian government aims to generate 27% of its energy from renewables by 2030, six grid-connected solar power plants have already been commissioned reaching a total capacity of 91 MW (Yaneva. M. 2016). The province of Adrar is also favorable for wind energy systems. Wind speed distributions show that 20% of wind speeds are 4 m/s or less, 62.5% are between 5~9 m/s, while the rest are above 10 m/s according to the wind atlas (Atlas Climatologique National. 2008). The latter states that average speeds are the lowest from 8 pm to 5 am. The maximum speeds are reached in the morning around 10 am and medium speeds are maintained throughout the day. The region of Adrar is characterized by a hot and very arid Saharan climate, cold nights in winter, scarce rainfall, and seasonal sandstorms. The hottest month is July where the mean minimum and maximum temperatures are 26.8 °C and 44.9 °C respectively. January marks the coldest month with 4.5 °C and 20.3 °C as the respective mean minimum and maximum temperatures (Atlas Climatologique National. 2008). Table 1 shows average horizontal global irradiation, average maximum temperatures and humidity levels on a monthly basis (Dabou et al. 2016).

**Table 1: Monthly average of daily weather parameters in 2014 (Dabou et al. 2016).**

| Months | H (kW h/m$^2$) | $T_{max}$ (C°) | RH (%) |
|---|---|---|---|
| January | 4.36 | 22.67 | 38 |
| February | 5.49 | 26.06 | 26 |
| March | 6.64 | 27.98 | 21 |
| April | 7.73 | 32.72 | 14 |
| May | 7.8 | 39.76 | 12 |
| Jun | 8.1 | 42.42 | 11 |
| July | 7.48 | 47.39 | 8 |
| August | 6.96 | 45.86 | 13 |
| September | 6.16 | 43.93 | 15 |
| October | 5.48 | 37.15 | 18 |
| November | 4.23 | 27.28 | 35 |
| December | 4.26 | 20.68 | 43 |
| Mean Monthly | 6.22 | 34.49 | 21.16 |

### 2.1.1. Weather data of Adrar

Data of solar irradiation, wind speed and temperature were extracted from a weather file generated using METEONORM 7.1. The data is in an hourly form and covers many parameters. The data covers an average of the ten last years up to 2021. Relevant data is converted into a monthly basis and represented in Fig. 4-6.

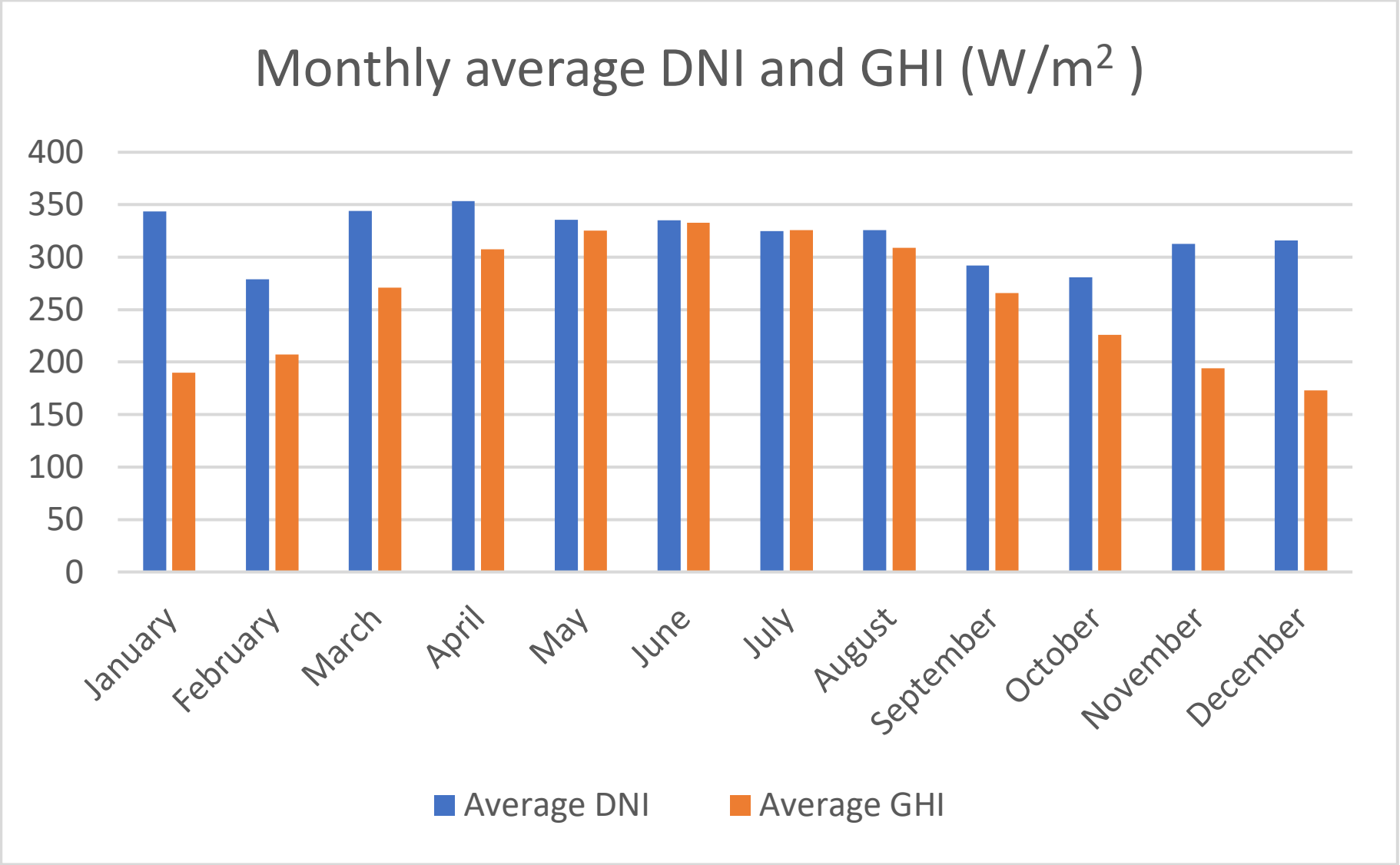

**Fig. 4: Monthly average direct normal and global horizontal irradiation.**

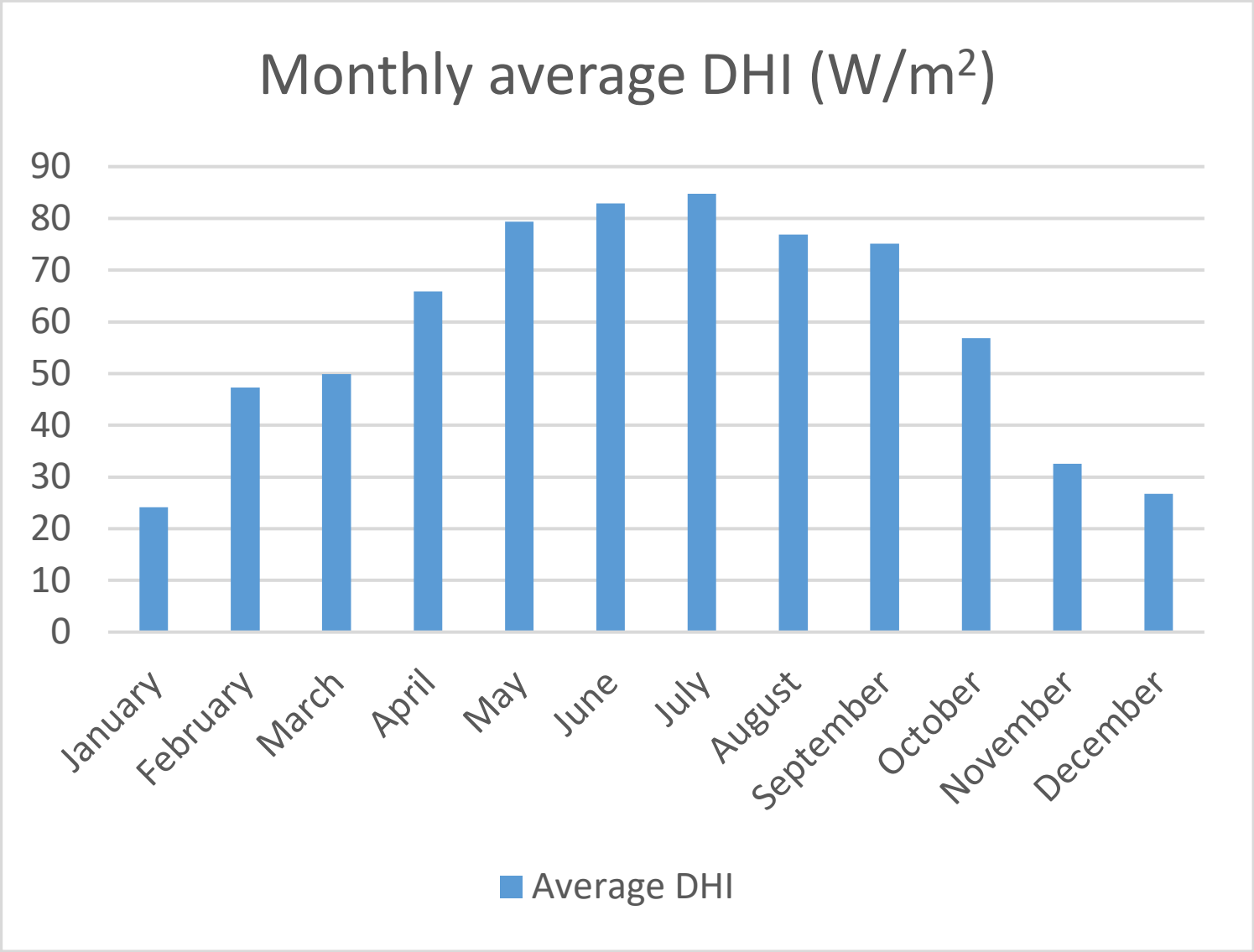

**Fig. 5: Monthly average diffuse horizontal irradiation.**

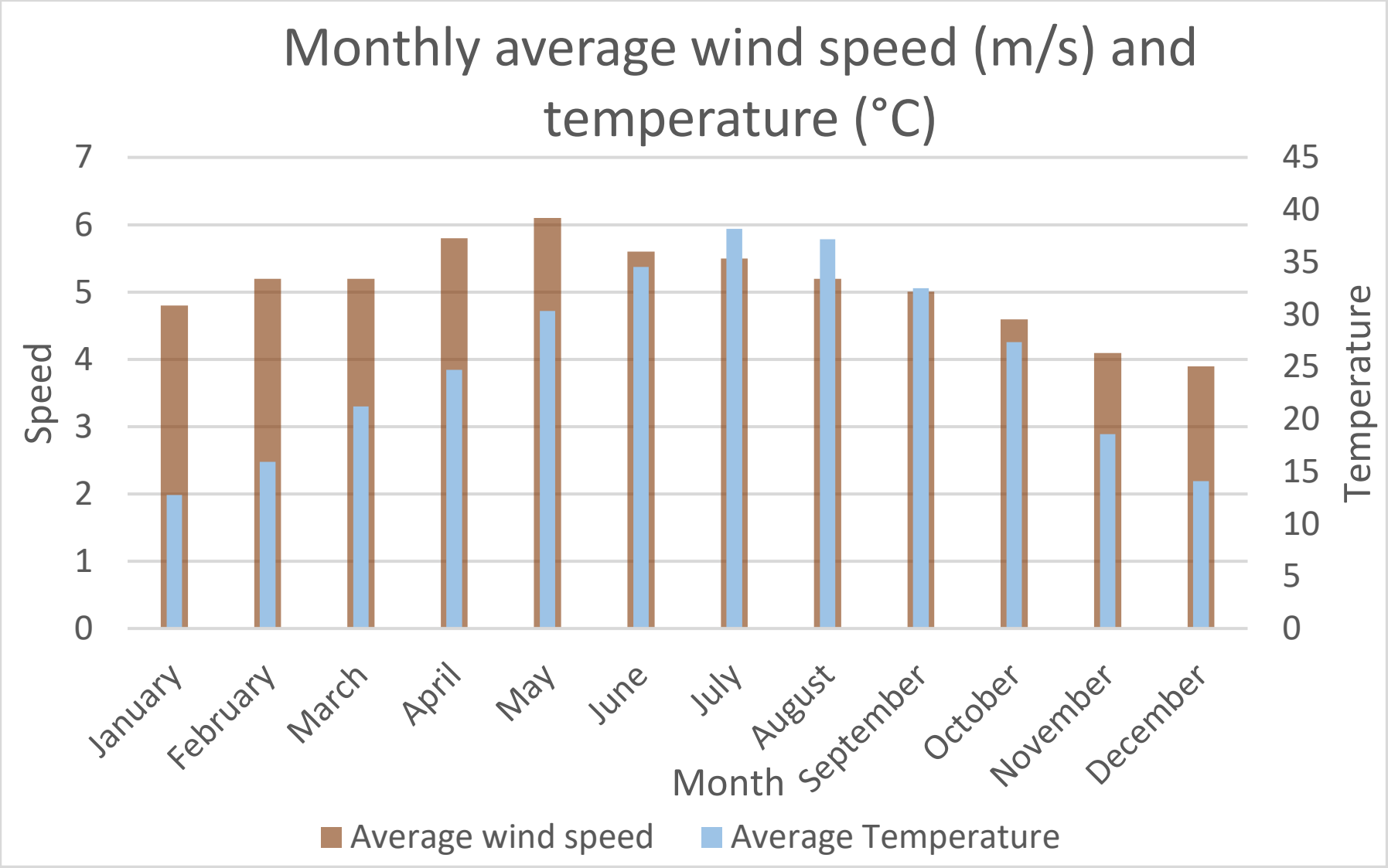

**Fig. 6: Monthly average wind speed and temperature.**

The figures show high average direct normal irradiation (DNI) content, especially in winter and spring, and the software indicates an error of 10%. The DNI levels tend to decrease in February and October. The monthly average direct horizontal irradiation (DHI) and global horizontal irradiation (GHI) levels follow a similar trend where the highest levels are recoded during the summer season, especially June and July. Average monthly wind speed reaches its maximum of around 6 m/s in May, and the maximum speed is around 16 m/s, which is very suitable for wind energy systems considering wind turbines with high hub heights. Overall, summer and spring seasons receive most of the wind energy potential. The trend for average temperature shows a noticeable increase during summer, the peak is recorded in July with and average value of 37 °C. This would decrease the efficiency of PV generators but for the other seasons, it is much cooler and more favorable. The software indicates a temperature error of 0.3 °C.

### 2.1.2. Load Profile of Adrar

The load profile of Adrar covers the amount of power consumed by the different sectors in the town. Daily data was used to determine the peak demand in the region. For representing the load, hourly power production data of the Adrar's central power plant was used instead, due to its availability. It is worth mentioning that the central power plant not only supplies the entire commune of Adrar, but part of its production goes to the neighboring cities as well. Fig. 7 depicts the monthly average power produced from the central power plant.

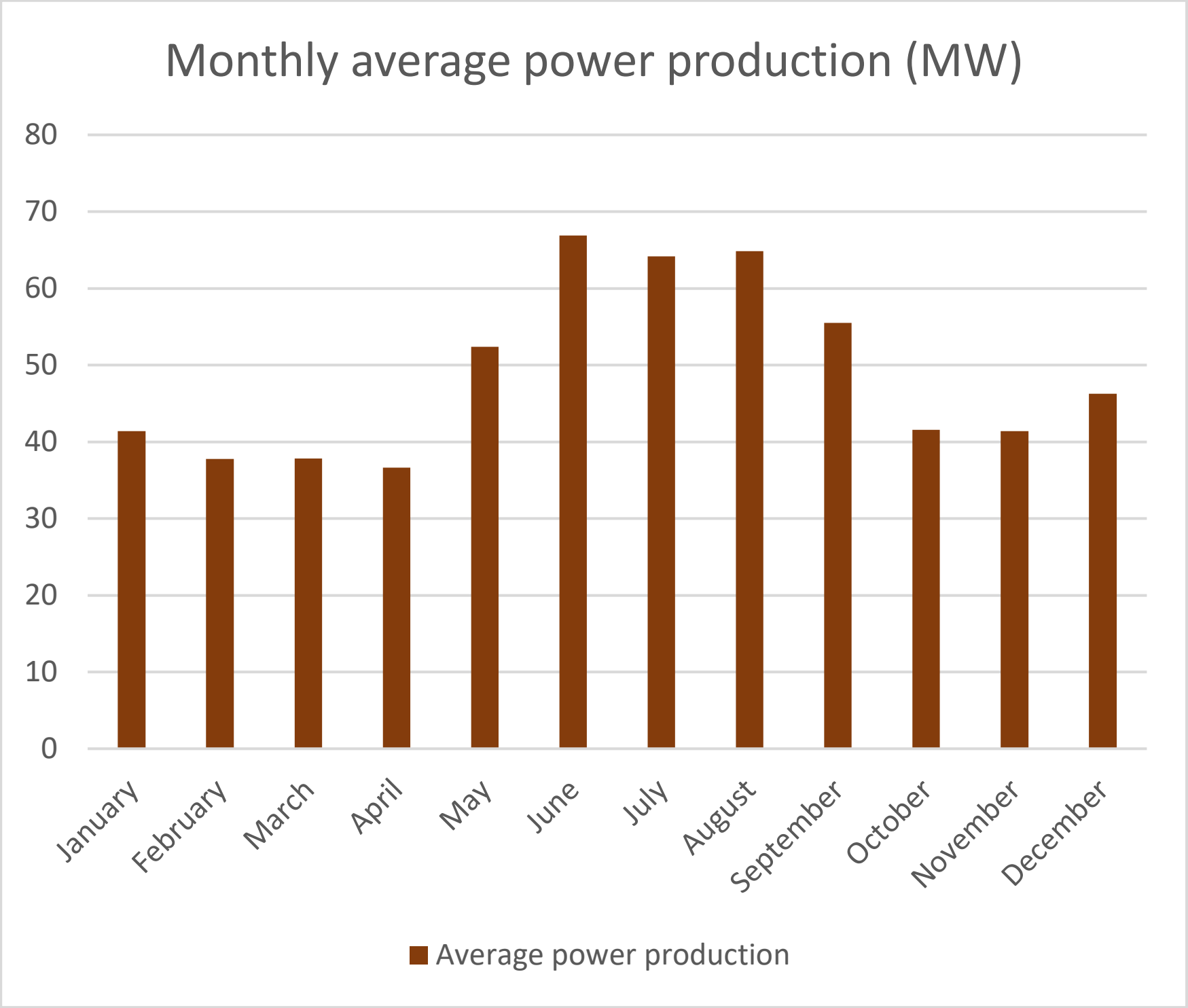

**Fig. 7: Monthly average power output of Adrar central power plant.**

The peak demand occurs during summer, this is justified by the cooling load since the place has a Saharan climate and temperatures skyrocket during this period. Although consumption rates in August are the highest, the peak demand for 2017 occurred in June reaching 63 MW, this represents an increase of about 6.5% compared to 2016. Based on the observed data, August and January are selected for simulation analysis in this study. August represents the summer period where the peak loads appear accompanied by high temperatures, while January represents the winter period where the wind energy potential and the load are the lowest. The choice of these two months aims to simulate the designed system for two extreme cases to be able to evaluate its contribution to power production.

## 2.2. System Design and Sizing

The sizing of the system is based on the production capacity and the energy potential of Adrar, while the design addresses the performance and capital cost of the overall system. Thus, the system is sized to meet at least half of the peak demand in Adrar, opting for a higher capacity is not viable since the capital cost will be too high and a vast area will be required. This has so many drawbacks, for instance, it will be too cumbersome to clean PV panels from dust depositions. In addition, the design employs certain topologies to maximize power output, minimize losses and reduce the upfront investment needed.

### 2.2.1. PV system

The PV generator produces 5 MW and consists of 14320 panels, with each panel having a rated power of 350 watts. To provide higher power levels, the multi-string inverter topology has been adopted. Each DC-DC boost converter is linked to an inverter, while connecting five PV-panel sets in parallel. Each PV-panel set contains eight panels in series. The whole PV system contains 358 inverters, 358 boost converters and a total of 14320 panels distributed over the 358 strings. To minimize the surface needed for panels, the two sets of each string face southwards with a spacing distance of three meters between the panels to avoid shading. The strings are then placed side by side as shown in Fig. 8, with this configuration the area is minimized to about 0.064 km². Based on the latitude of Adrar ($25° < \emptyset=27.52° < 50°$) and according to (SolarPanelTilt.com 2021), the optimal tilt angle for the PV panels is obtained as: ($\beta= 0.76 \times \emptyset + 3.1=25°$).

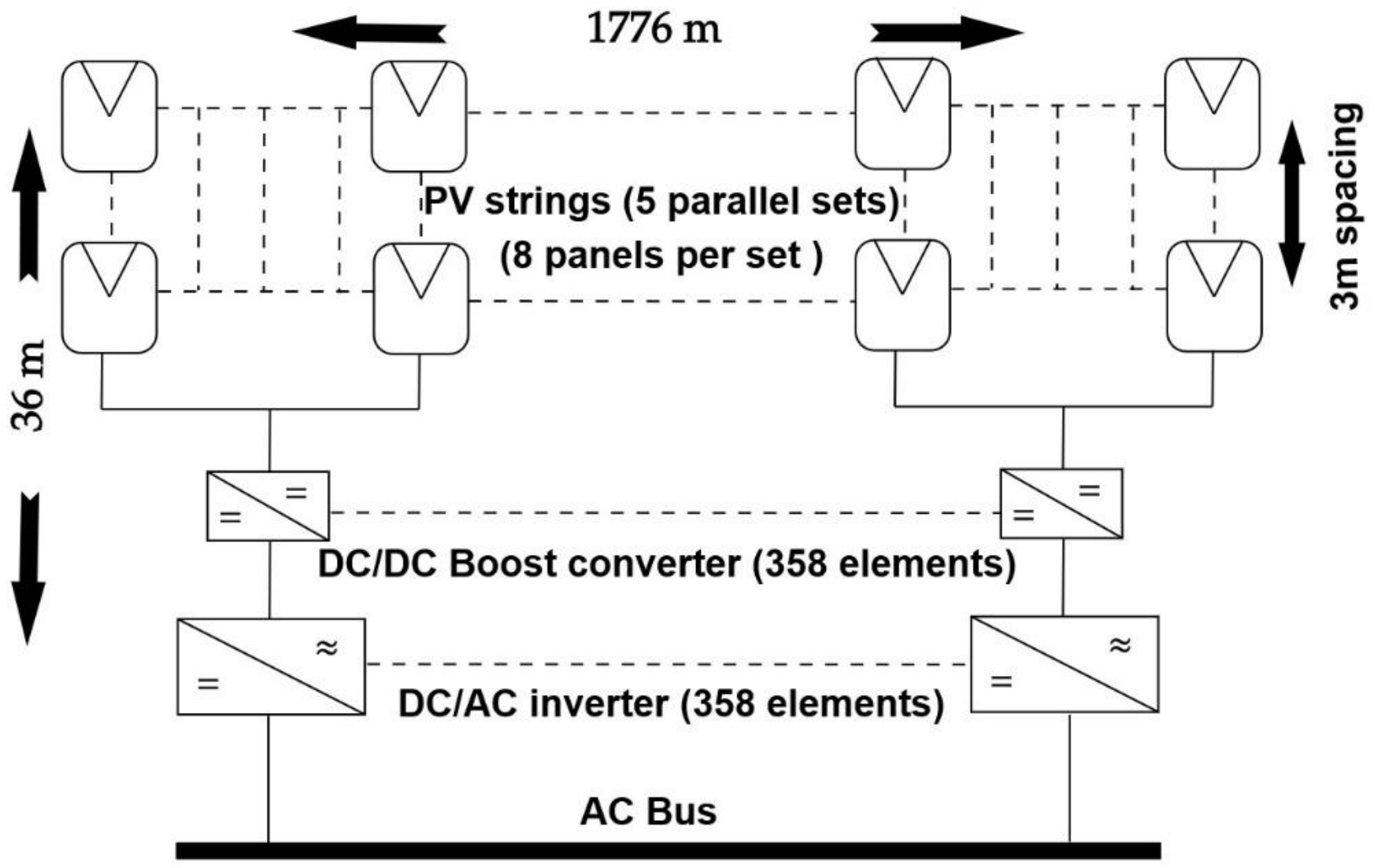

**Fig. 8: PV plant scheme**

### 2.2.2. Wind system

The wind farm consists of 10 large wind turbines; each one has a rated power of 3 MW. Adrar is a region of pressure difference close to the Atlantic Ocean where the wind blows from the west. Thus, the wind turbines face westwards. The general rule-of-thumb for wind farm spacing speculates a spacing distance of about five to seven rotor diameters between two turbines (Meyers et al. 2012). Hence, to avoid turbulence and ensure maximum power extraction, a spacing distance of 450 meters is maintained between the turbines. The 10 turbines are divided into three sets. The first set comprises four aligned turbines, while the other two sets contain three turbines each as shown in Fig. 9. The estimated surface required for the wind farm is about 1.8 km².

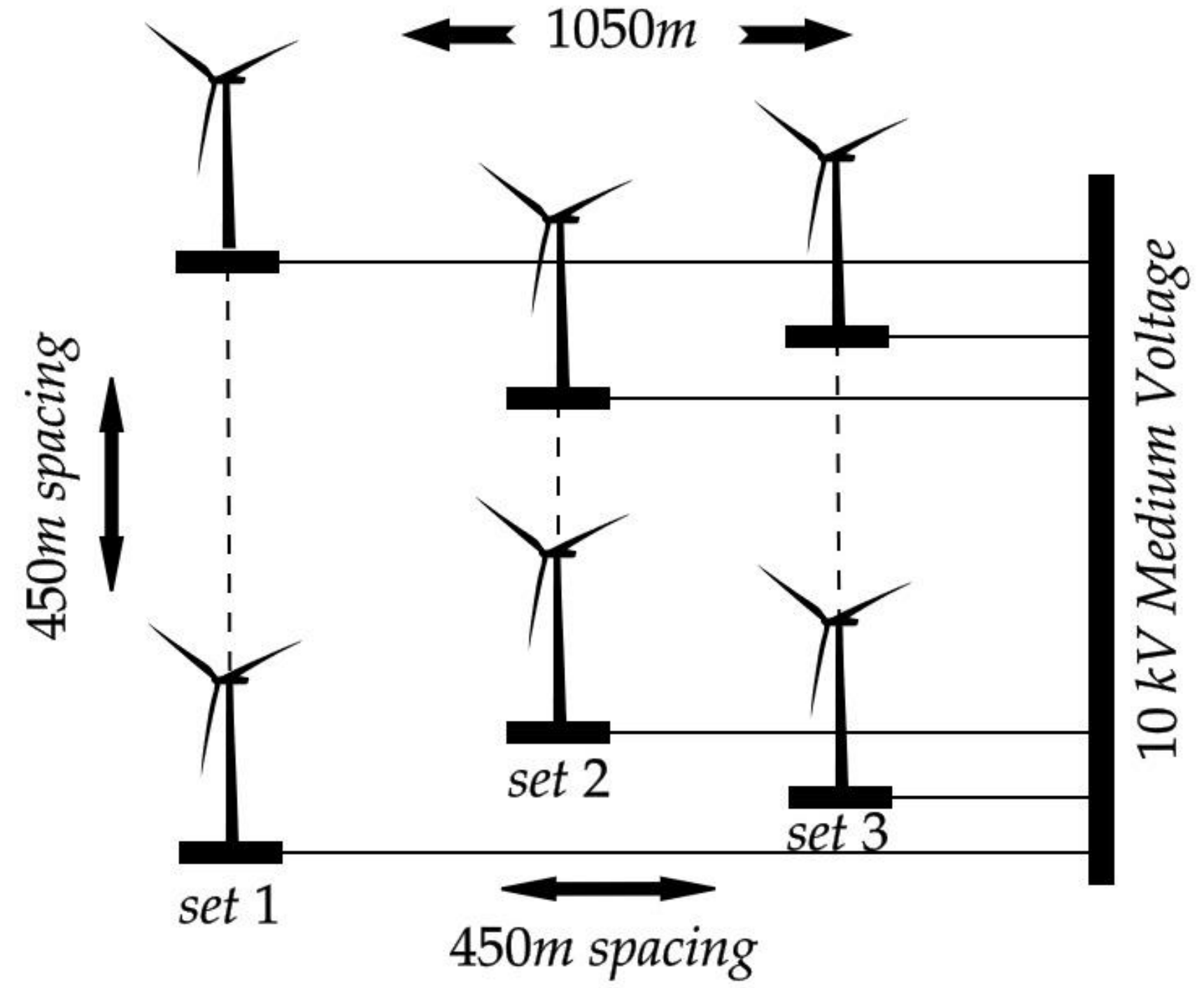

**Fig. 9: Wind farm scheme**

The wind generators are linked to medium voltage lines (10 kV), a step-up transformer is required to match the output voltage with that of the medium voltage lines and allow grid connection. The total system has a rated power of 35 MW, but this is rarely maintained due to the fluctuations of irradiation levels and wind speeds. However, the system should undoubtedly improve the grid performance and eliminate power outages during peak demands.

## 2.3. Modelling of the PV system

### 2.3.1. PV generator

The PV array is modeled based on four equations that govern the output current [68]:

Photo-current: $$I_{ph} = [I_{sc} + k_i.(T - T_n)].\frac{G}{1000} \quad (1)$$

Saturation current: $$I_0 = I_{rs}.\left(\frac{T}{T_n}\right)^3.\exp\left[\frac{q.E_{g0}.\left(\frac{1}{T_n}-\frac{1}{T}\right)}{n.K}\right] \quad (2)$$

Reverse saturation current: $$I_{rs} = \frac{I_{sc}}{\exp\left(\frac{q.V_{oc}}{n.N_s.K.T}\right)-1} \quad (3)$$

Output current: $$I = I_{ph} - I_0.\left[\exp\left(q.\frac{(V+I.R_s)}{n.K.N_s.T}\right) - 1\right] - \frac{V+R_s.I}{R_p} \quad (4)$$

Where $I_{ph}$, $I_0$, $I_{rs}$ and $I_{sc}$ are the respective photogenerated current, diode saturation current, reverse saturation current and short-circuit current. $k_i$, $T$, $T_n$, $q$ and $G$ are the current coefficient, the temperature, the nominal temperature, the electron charge and the irradiation respectively. $E_{g0}$, $K$, $N_s$, $R_s$ and $R_p$ represent the semiconductor bandgap energy, the Boltzmann's constant, the number of cells connected in series and the equivalent series and parallel resistors respectively. The monocrystalline Amerisolar AS-6M PV panel was chosen due to its high efficiency of 18.04%, competitive price of 0.3 $\kW and good performance in hot climates. Moreover, monocrystalline solar cells proved better performance in desert climates than other types of solar cells according to (Mosalam et al. 2000). The main panel specifications are given in Table 2 (see Appendix). The series ($R_s$) and shunt ($R_p$) resistances are calculated according to (Villalva et al. 2009).

### 2.3.2. Maximum power point control

To keep the PV generator working at maximum power, an MPPT tracker has been added to the model. Perturb and observe algorithm has been chosen, because of its precision, simplicity and ease of implementation. The algorithm is based on the introduction of perturbations to the voltage and the observation of power output change. The signs of the last perturbation and the last power increment are used to decide the next perturbation; this continuous loop keeps the output power at its maximum. The output voltage of the array is injected directly into a DC-DC converter. Fig. 10 shows the MPPT algorithm.

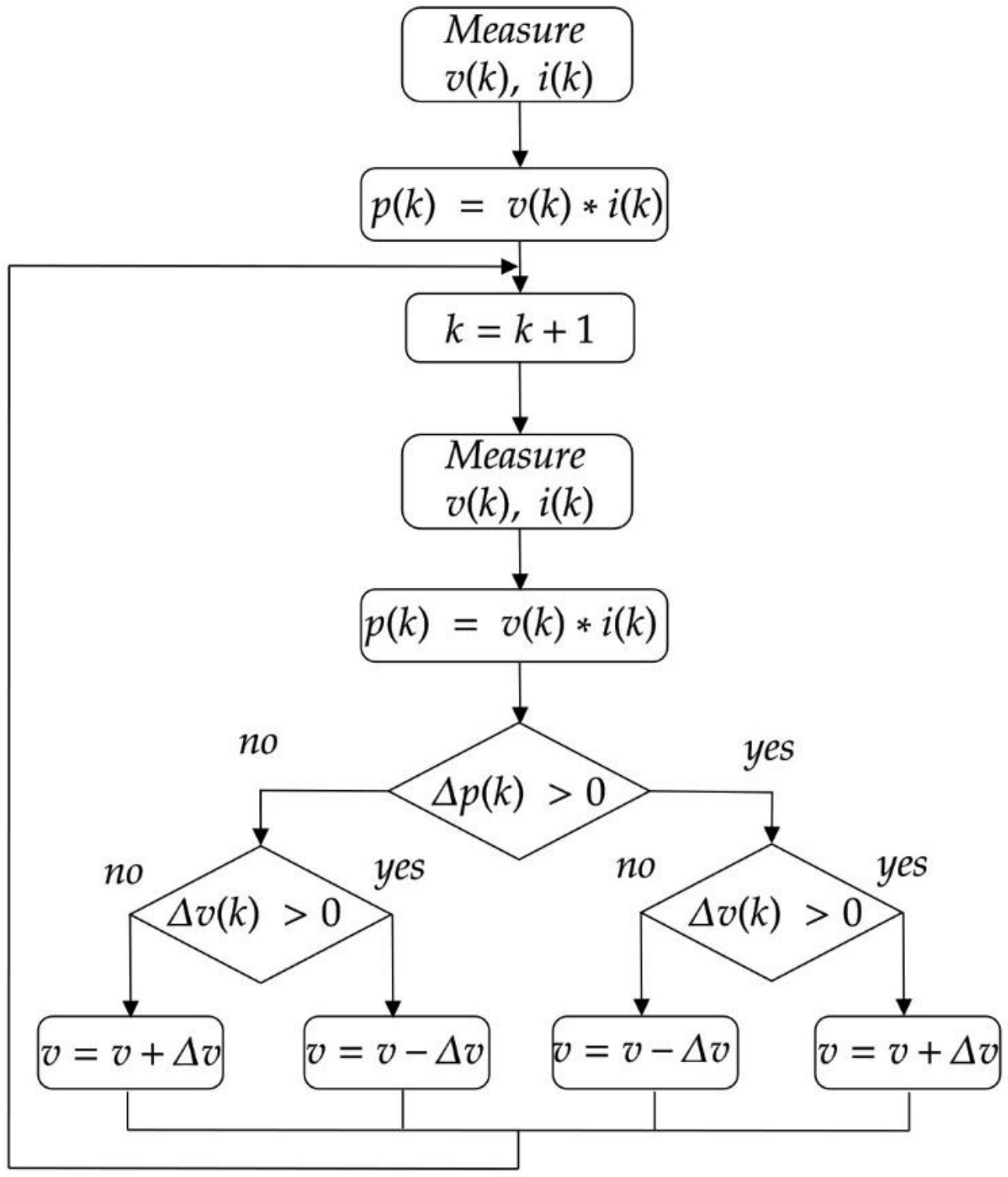

**Fig. 10 MPPT algorithm**

### 2.3.3. DC-DC converter

The DC boost converter has high efficiency and robustness, and it is capable of stepping up the output. The converter receives the fluctuating DC output voltage of the PV array and converts it into a fixed stepped-up DC voltage. The equivalent electrical circuit is shown in Fig. 11, and equations (5,6) are used for modeling.

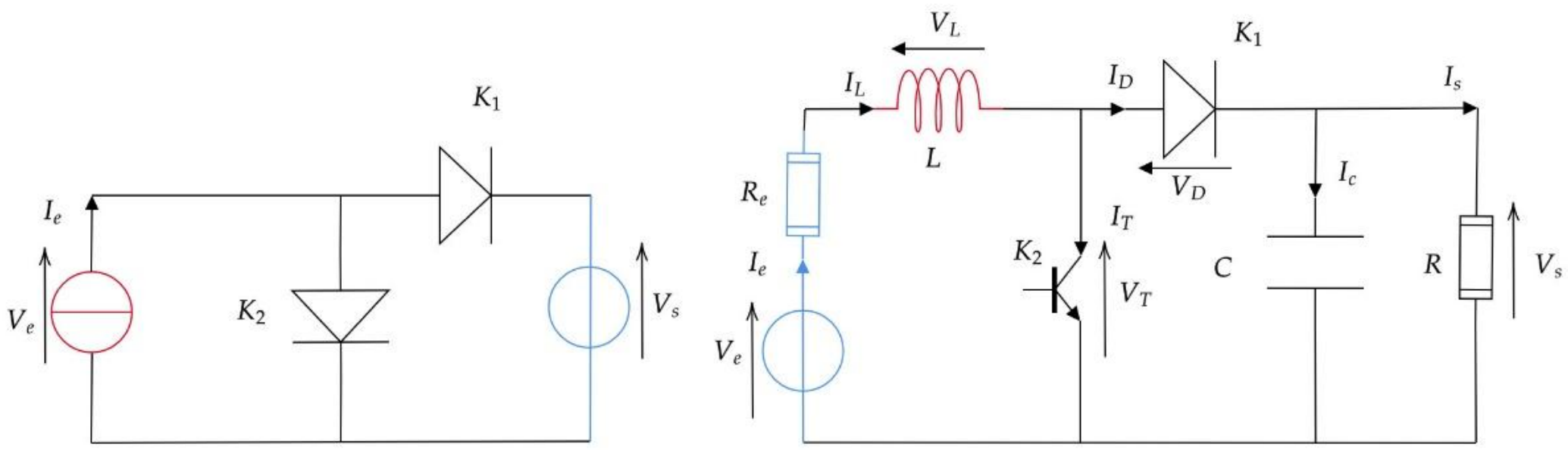

**Fig. 11: DC-DC boost converter.**

Loop 1: $$L\frac{di_L}{dt} = V_e - V_c(1-\alpha) \quad (5)$$

Loop 2: $$C\frac{dV_c}{dt} = i_L(1-\alpha) - \frac{V_c}{R} \quad (6)$$

Where $\alpha$ represents the duty cycle obtained from the MPPT and used to generate the Pulse width modulation (PWM) signal which then controls the opening and closing of the switch (represented by $K_2$ in the circuit).

### 2.3.4. Three-phase inverter

The connection of the PV system to the AC bus is done via a voltage source inverter. The output of the boost converter is a fixed DC voltage; it is fed into a grid-tied three-phase inverter. Such inverters are used only in high power applications and their output is a three-phase voltage ($v_a$, $v_b$ and $v_c$), the latter is controlled by PWM signals. Pulse Width modulation is the process of changing the width of a pulse train in direct proportion to a small control signal. The electric circuit of the inverter with RL filter is shown in Fig. 12.

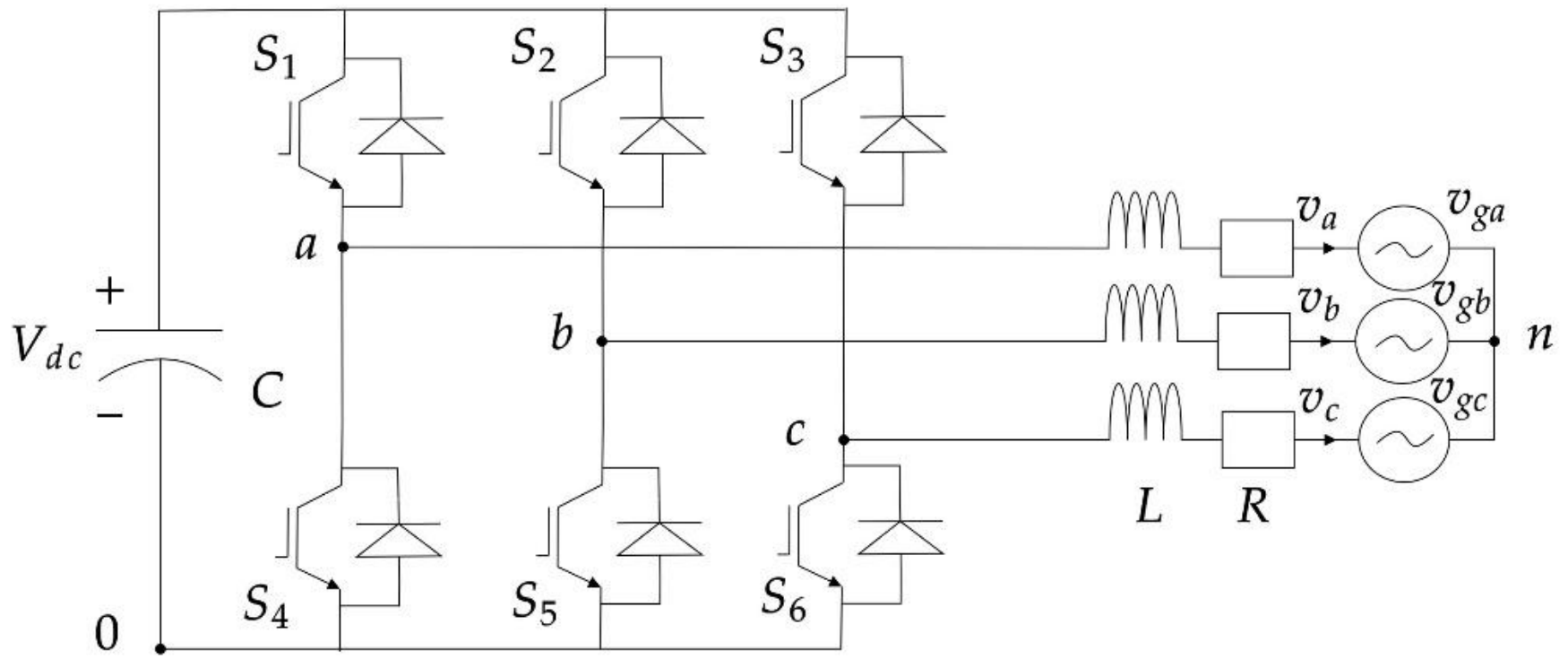

**Fig. 12: Equivalent electric circuit of the three-phase inverter.**

Three-phase inverters can be controlled in three types of frames, the synchronous rotating frame (d-q), the stationary reference frame ($\alpha$-$\beta$) and the natural frame (a-b-c). In d-q frame control, three-phase voltages and currents are transformed by Park's transformation into the d-q reference frame that rotates synchronously with the grid voltage. Thus, three-phase variables become DC quantities. As the control variables are DC, different filtering methods can be used. In addition, an innovative PSO-PI, which will be detailed, is designed to achieve good performance in this reference frame control. The model of the inverter is determined by the equations (7-9) (Isen et al. 2016):

a-b-c $$V_{k0} = L\frac{di_k}{dt} + Ri_k + V_{gk} + V_{n0} \quad (7)$$

$$V_{n0} = \frac{1}{3}(V_{a0} + V_{b0} + V_{c0}) \quad (8)$$

$$V_{k0} = S_k\, V_{dc} = \begin{cases} V_{dc} & S_k = 1 \\ 0 & S_k = 0 \end{cases} \quad (9)$$

Where k, $S$ and R represent the phase of the inverter, the state of upper switches and the filter series resistor respectively. L and $V_n$ are the filter inductance and the voltage difference between grid neutral point and negative $V_{dc}$. In each phase of the inverter, the two switches are supposed perfect and work simultaneously, this allows us to write equation (10). Considering the voltages between a, b, c and the ground, we can express the model of the inverter by equation (11) (Loucif. M. 2016):

$$\begin{bmatrix} v_{a0} \\ v_{b0} \\ v_{c0} \end{bmatrix} = \frac{V_{dc}}{2} \times \begin{bmatrix} S_a \\ S_b \\ S_c \end{bmatrix} \quad (10)$$

$$\begin{bmatrix} v_a \\ v_b \\ v_c \end{bmatrix} = \frac{1}{3} \times \begin{bmatrix} 2 & -1 & -1 \\ -1 & 2 & -1 \\ -1 & -1 & 2 \end{bmatrix} \times \frac{V_{dc}}{2} \times \begin{bmatrix} S_a \\ S_b \\ S_c \end{bmatrix} \quad (11)$$

To control the inverter, the grid voltages are measured and used to determine the grid angle θ by means of Phase Lock Loop (PLL). (Isen et al. 2016) explains the PLL method in detail. The grid angle is then used in Park's transformation to transform the grid voltages ($V_{ga}$, $V_{gb}$ and $V_{gc}$) from a-b-c reference frame to d-q frame which simplifies the control, the result is $V_{gd}$ and $V_{gq}$ (the d-q components of the grid voltages). Following the same transformation method, $V_{id}$ and $V_{iq}$ (the d-q inverter components) are generated from the output voltages of the inverter. Park's transformation is given as (Isen et al. 2016):

$$\begin{bmatrix} V_d \\ V_q \end{bmatrix} = \frac{2}{3} \begin{bmatrix} \cos\theta & \cos\left(\theta - \frac{2\pi}{3}\right) & \cos\left(\theta - \frac{4\pi}{3}\right) \\ -\sin\theta & -\sin\left(\theta - \frac{2\pi}{3}\right) & -\sin\left(\theta - \frac{4\pi}{3}\right) \end{bmatrix} \begin{bmatrix} V_a \\ V_b \\ V_c \end{bmatrix} \quad (12)$$

### 2.3.5. RL Filter

The d-q voltage components of both the inverter and the grid are used in a vector control mechanism implemented within an RL filter to generate active (P) and reactive (Q) powers. The RL filter model is based on the following equations (Gaillard. A. 2010):

a-b-c:

$$V_a = -R_f.i_a - L_f\frac{di_a}{dt} + V_{ga} \quad (13)$$

$$V_b = -R_f.i_b - L_f\frac{di_b}{dt} + V_{gb} \quad (14)$$

$$V_c = -R_f.i_c - L_f\frac{di_c}{dt} + V_{gc} \quad (15)$$

d-q

$$V_d = -R_f.i_d - L_f\frac{di_d}{dt} + L_f\omega i_q + V_{gd} \quad (16)$$

$$V_q = -R_f.i_q - L_f\frac{di_q}{dt} - L_f\omega i_d + V_{gq} \quad (17)$$

Where:

- $L_f$ represents the inductance of the filter.
- $R_f$ is the resistor of the filter.

The RL filter connects the inverter to the grid because they have different voltage amplitudes and different phases, the filter maintains a voltage difference between the grid and the inverter which allows power to flow towards the grid. In addition, it reduces the harmonics that are generated by the electronic components of the system. The measured active and reactive powers are then directly compared to reference values (P*) and (Q*),

Q* is always set to zero to eliminate reactive power generation, this ensures zero phase angle between voltage and current leading to unity power factor. On the other hand, P* is set based on the output of the PV generator. The resulting errors are regulated by the PSO-PI controller and yield reference d-q inverter voltages ($V_d$* and $V_q$*). The control mechanism responsible for finding ($V_d$* and $V_q$*) is shown in Fig. 13:

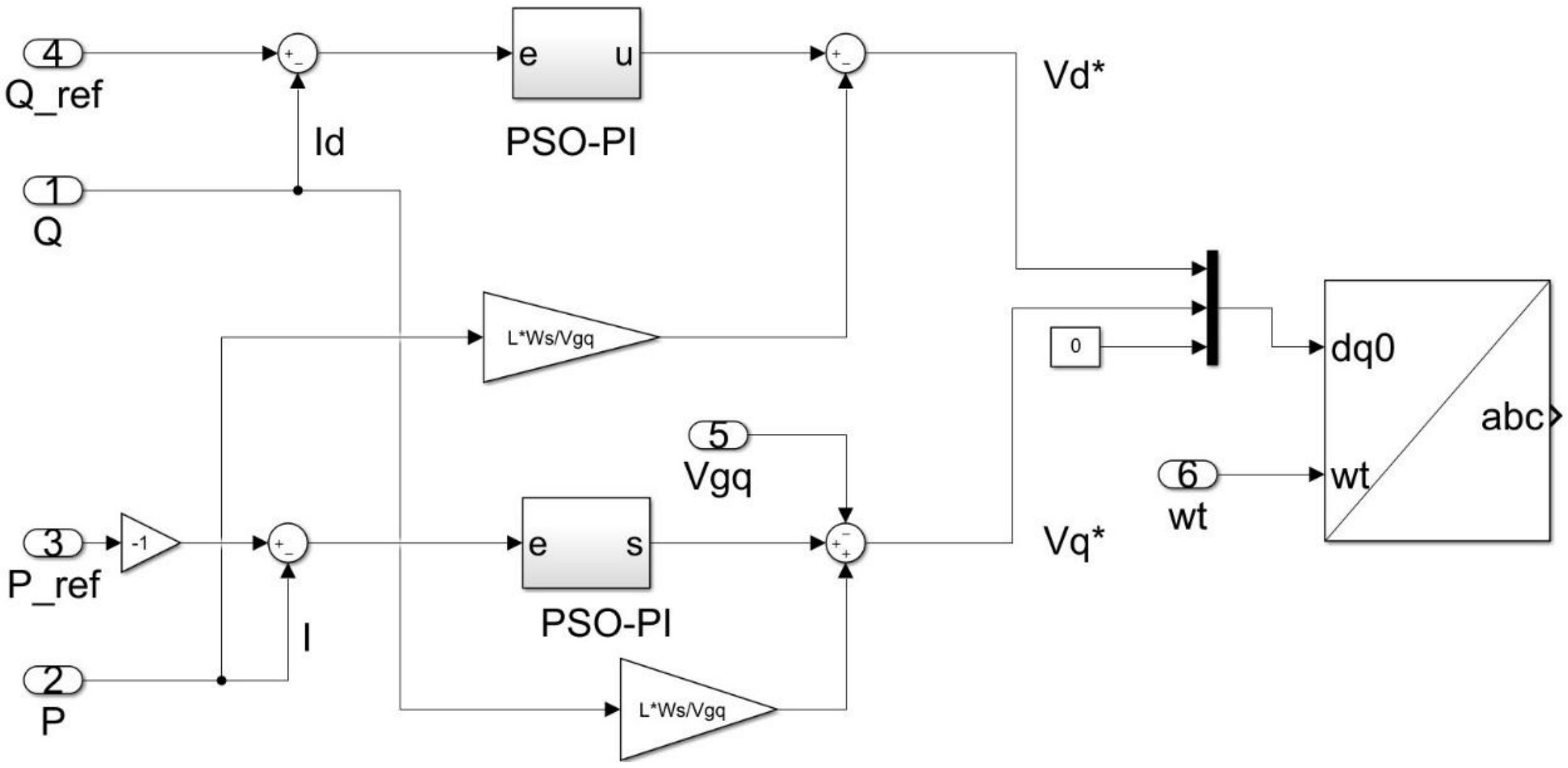

**Fig. 13: Simulink model of the inverter control**

With the help of Park's inverse transformation, the obtained voltages are transformed back into three-phase reference inverter voltages ($V_a$*, $V_b$* and $V_c$*) to be compared with a repetitive MATLAB sequence which in turn generates PWM signals that control the switches of the inverter. The PWM signals represent the state of the upper switches given in equation (9). Park's inverse transformation, known as d-q to a-b-c transformation is given by the following equation (Isen et al. 2016):

$$\begin{bmatrix} V_a \\ V_b \\ V_c \end{bmatrix} = \begin{bmatrix} \cos(\theta) & -\sin(\theta) \\ \cos\left(\theta - \frac{2\pi}{3}\right) & -\sin\left(\theta - \frac{2\pi}{3}\right) \\ \cos\left(\theta - \frac{4\pi}{3}\right) & -\sin\left(\theta - \frac{4\pi}{3}\right) \end{bmatrix} \begin{bmatrix} V_d \\ V_q \end{bmatrix} \tag{18}$$

### 2.4. Modelling of the wind system

The wind farm consists of a set of wind turbines, which are divided into dynamic and mechanical models, in addition to the electrical model of a doubly-fed induction generator (DFIG) linked to the turbine via a gearbox. The stator of the DFIG is directly connected to and synchronized with the grid, while the connection of the rotor is ensured by a back-to-back converter consisting of a rectifier, an inverter, a DC bus, and an RL filter as seen in Fig. 14.

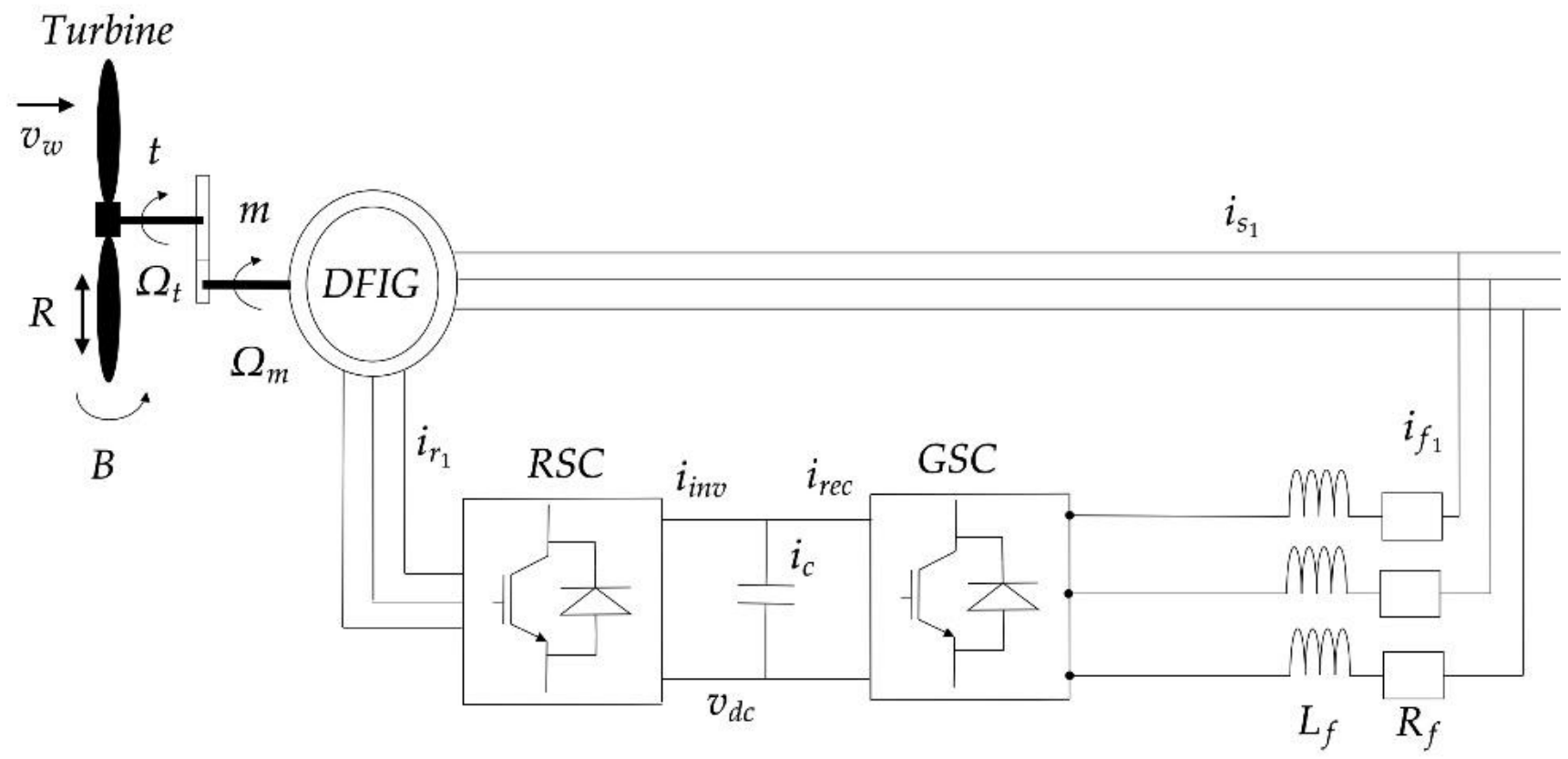

**Fig. 14: DFIG-based wind system.**

### 2.4.1. Dynamic model of the wind turbine

The mechanical power extracted from wind is characterized by equation (19):

$$P_t = \frac{1}{2}\rho Cp(\lambda,\beta)\pi R^2 V^3 \qquad (19)$$

with the TSR defined as:

$$\lambda = \frac{\Omega_t . R}{V} \qquad (20)$$

The parameters $\rho$, $Cp$ and $R$ are the air density, the power coefficient and the turbine blade radius respectively. $V$, $\Omega_t$, and $\beta$ are the respective wind velocity, rotational speed and pitch angle of the turbine. The maximum power coefficient is a function of λ and β. According to (Gaillard. A. 2010), the maximum power coefficient for a 3 MW wind turbine it is given by:

$$C_p = C_1 \times \left(C_2 \times \frac{1}{\lambda_i} - C_3 \times \beta - C_4\right) \times \exp\left(-\frac{C_5}{\lambda_i}\right) + C_6 \times \lambda \qquad (21)$$

Where ($\frac{1}{\lambda_i} = \frac{1}{\lambda + 0.08\times\beta} - \frac{0.035}{\beta^3+1}$) and the coefficients $C_i$ are given in Table 3 (see Appendix)

### 2.4.2. Mechanical model of the turbine

The torque of the generator ($T_m$) and turbine ($T_t$) are used to evaluate the speed ($\Omega_m$). The behavior of the mechanical system is governed by the classical rotational dynamics equation:

$$\left(\frac{J_t}{G^2} + J_m\right)\frac{d\Omega_m}{dt} + f_v \times \Omega_m = T_t - T_m \qquad (22)$$

Where $J_t$ and $J_m$ are the inertia moments of the turbine and the generator, and $f_v$ is the friction coefficient of the generator.

### 2.4.3. MPPT control for the wind turbine

The MPPT technique relies on controlling the rotational speed of the generator to ensure maximum power extraction, this is achieved by imposing a reference torque on the machine. Thus, the TSR must be maintained at its optimal value ($\lambda_{opt}$) to maintain the power coefficient at its maximum value $C_{pmax}$ as well. The turbine control and the MPPT technique are demonstrated in Fig. 15.

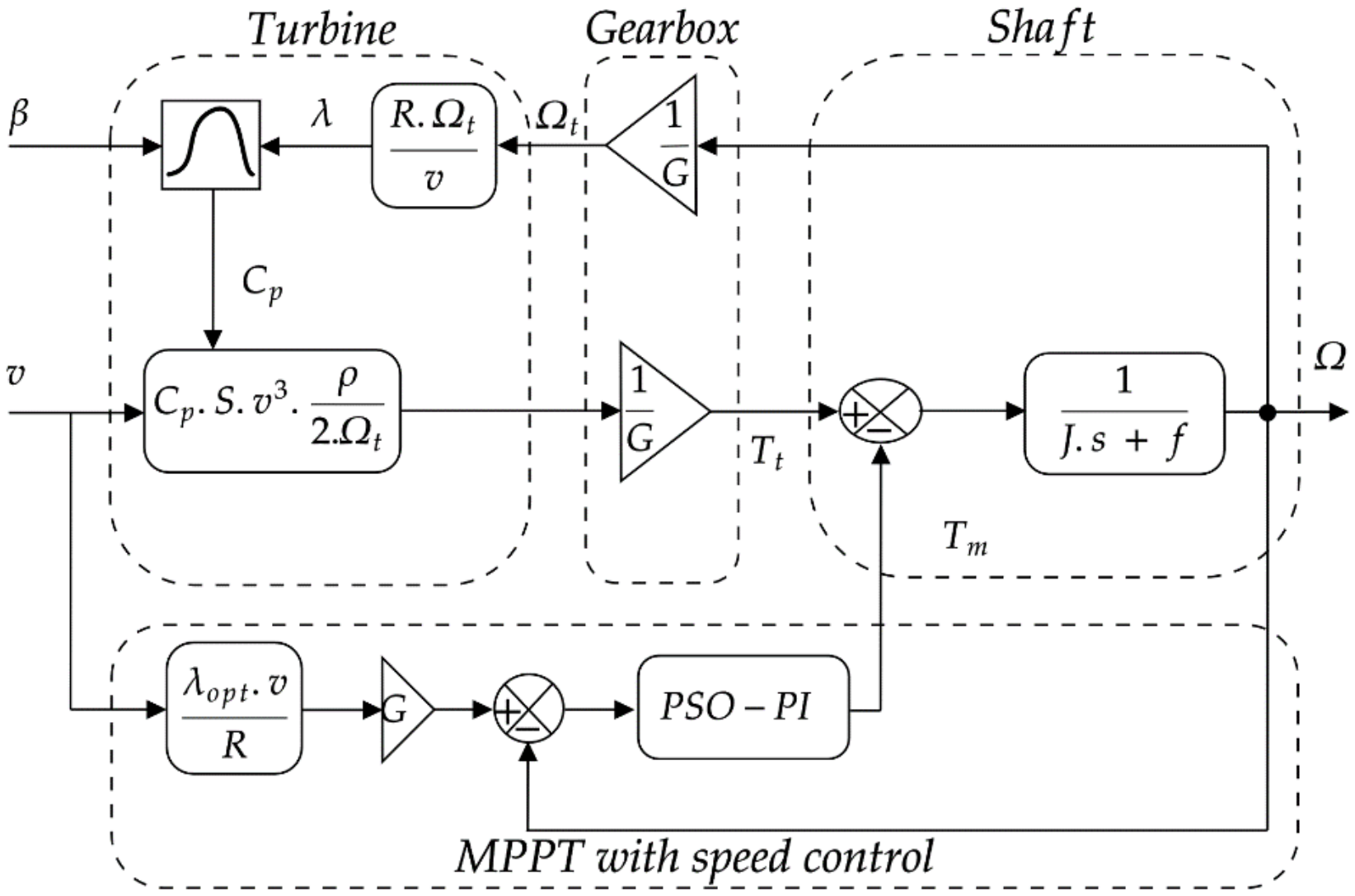

**Fig. 15: Wind turbine with MPPT control mechanism.**

### 2.4.4. Model of the doubly-fed induction generator

The DFIG is modeled in Park's reference frame to facilitate the implementation of the vector control method. The model described in (Gaillard. A. 2010) is simplified where the fluxes are additive and the inductances are constant with a mutual sinusoidal variation between the rotor and stator windings as a function of the electrical angle $\theta_e$ and their axes. Equations (23-26) govern the stator and rotor voltages in the Park reference frame.

$$V_{sd} = R_s . i_{sd} + \frac{d\varphi_{sd}}{dt} - \dot{\theta_s} . \varphi_{sq} \quad (23)$$

$$V_{sq} = R_s . i_{sq} + \frac{d\varphi_{sq}}{dt} + \dot{\theta_s} . \varphi_{sd} \quad (24)$$

$$V_{rd} = R_r . i_{rd} + \frac{d\varphi_{rd}}{dt} - \dot{\theta_r} . \varphi_{rq} \quad (25)$$

$$V_{rq} = R_r . i_{rq} + \frac{d\varphi_{rq}}{dt} - \dot{\theta_r} . \varphi_{rd} \quad (26)$$

Where $V_{sd}$ ,$V_{sq}$ and $V_{rd}$, $V_{rq}$ are the stator and rotor d-q voltages, $i_{sd}$, $i_{sq}$ and $i_{rd}$, $i_{rq}$ are the stator and rotor d-q currents. $\varphi_{sd}$, $\varphi_{sq}$ and $\varphi_{rd}$, $\varphi_{rq}$ are the stator and rotor d-q fluxes, $R_r$ and $R_s$ are the rotor and stator resistances and $\theta_r$ and $\theta_s$ are rotor and stator angles. The stator and rotor flux d-q components are given by:

$$\varphi_{sd} = L_s . i_{sd} + m. L_m . i_{rd} \quad (27)$$

$$\varphi_{sq} = L_s . i_{sq} + m. L_m . i_{rq} \quad (28)$$

$$\varphi_{rd} = L_r . i_{rd} + m. L_m . i_{sd} \quad (29)$$

$$\varphi_{rq} = L_r . i_{rq} + m. L_m . i_{sq} \quad (30)$$

$$L_s = L_{fs} + L_m \quad (31)$$

$$L_r = L_{fr} + m^2 . L_m \quad (32)$$

Where $L_s$, $L_r$ and $L_m$ are the stator/rotor cyclic inductances and the magnetizing inductance respectively. $L_{fs}$ and $L_{fr}$ are the stator and rotor leakage inductances and $m$ is the transformation ratio. Formulas (34-37) are used to calculate the stator and rotor active and reactive powers.

$$P_S = V_{sd}.i_{sd} + V_{sq}.i_{sq} \quad (34)$$
$$Q_s = V_{sq}.i_{sd} - V_{sd}.i_{sq} \quad (35)$$
$$P_r = V_{rd}.i_{rd} + V_{rq}.i_{rq} \quad (36)$$
$$Q_r = V_{rq}.i_{rd} - V_{rd}.i_{rq} \quad (37)$$

The electromagnetic torque $T_m$ can be expressed in equation (38-39):

$$T_m = p.(\varphi_{sd}.i_{sq} - \varphi_{sq}.i_{sd}) \quad (38)$$
$$T_m = p.\frac{m.L_m}{L_s}.(\varphi_{sq}.i_{rd} - \varphi_{sd}.i_{rq}) \quad (39)$$

Where p represents the number of pole pairs.

### 2.4.5. DFIG Control

By neglecting the resistance of the stator windings, the simplified model of the DFIG in the d-q reference frame is obtained as the following:

$$V_{sd} = 0 \quad (40)$$
$$V_{sq}istance = U_s = \omega_s.\varphi_{sd} \quad (41)$$
$$V_{rd} = R_r.i_{rd} + \frac{d\varphi_{rd}}{dt} - \omega_r.\varphi_{rq} \quad (42)$$
$$V_{rq} = R_r.i_{rq} + \frac{d\varphi_{rq}}{dt} - \omega_r.\varphi_{rd} \quad (43)$$

The stator currents are obtained from the stator and rotor flux equations as:

$$i_{sd} = \frac{\varphi_{sd} - m.L_m.i_{rd}}{L_s} \quad (44)$$
$$i_{sq} = -m.\frac{L_m}{L_s}.i_{rq} \quad (45)$$

Thus, rotor flux equations become like the following:

$$\varphi_{rd} = \left(L_r - \frac{(m.L_m)^2}{L_s}\right).i_{rd} + m.\frac{L_m}{L_s}.\varphi_{sd} = \sigma.L_r.i_{rd} + m.\frac{L_m}{L_s}.\varphi_{sd} \quad (46)$$
$$\varphi_{rq} = L_r.i_{rq} - \frac{(m.L_m)^2}{L_s}.i_{rq} = \sigma.L_r.i_{rq} \quad (47)$$

With $\sigma$ being the leakage coefficient:

$$\sigma = 1 - \frac{(m.L_m)^2}{L_s.L_r} \quad (48)$$

By replacing the direct and quadrature components of the rotor fluxes we obtain:

$$v_{rd} = R_r.i_{rd} + \sigma.L_r\frac{di_{rd}}{dt} + e_{rd} \quad (49)$$
$$v_{rq} = R_r.i_{rq} + \sigma.L_r\frac{di_{rq}}{dt} + e_{rq} + e_\varphi \quad (50)$$

Where:

$$e_{rd} = -\sigma.L_r.\omega_r.i_{rq} \quad (51)$$
$$e_{rq} = \sigma.L_r.\omega_r.i_{rd} \quad (52)$$
$$e_\varphi = \omega_r.m.\frac{L_m}{L_s}.\varphi_{sd} \quad (53)$$

The electromagnetic torque becomes as follows:

$$T_m = -p.\frac{m.L_m}{L_s}.\varphi_{sd}.i_{rq} \quad (54)$$

Hence, the active and reactive powers are given by equations (55,56):

$$P_S = -v_{sq}.m.\frac{L_m}{L_s}.i_{rq} \quad (55)$$

$$Q_S = \frac{v_{sq}.\varphi_{sd}}{L_s} - v_{sq}.m.\frac{L_m}{L_s}.i_{rd} \quad (56)$$

From these expressions it is obvious that the choice of the d-q frame makes power produced from the stator proportional to the rotor q current component, while the reactive power is not proportional to the rotor d current component. This implies that independent control of the stator's active and reactive power is possible through the rotor d-q current components. For controlling the rotor d-q current components, reference d-q rotor currents are needed. These are obtained from the stator flux, which is estimated by the following formula:

$$\varphi_{sd-est} = L_s.i_{sd} + m.L_m.i_{rd} \quad (57)$$

The reference rotor current components are then generated as follows:

$$i_{rq}^* = -\frac{L_s}{p.m.L_m.\varphi_{sd-est}}.T_m^* \quad (58)$$

$$i_{rd}^* = \frac{\varphi_{sd-est}}{m.L_m} - \frac{L_s}{m.L_m.v_{sq}}.Q_s^* \quad (59)$$

Fig. 16 shows the control mechanism of the DFIG, ensured by the rotor side converter (RSC).

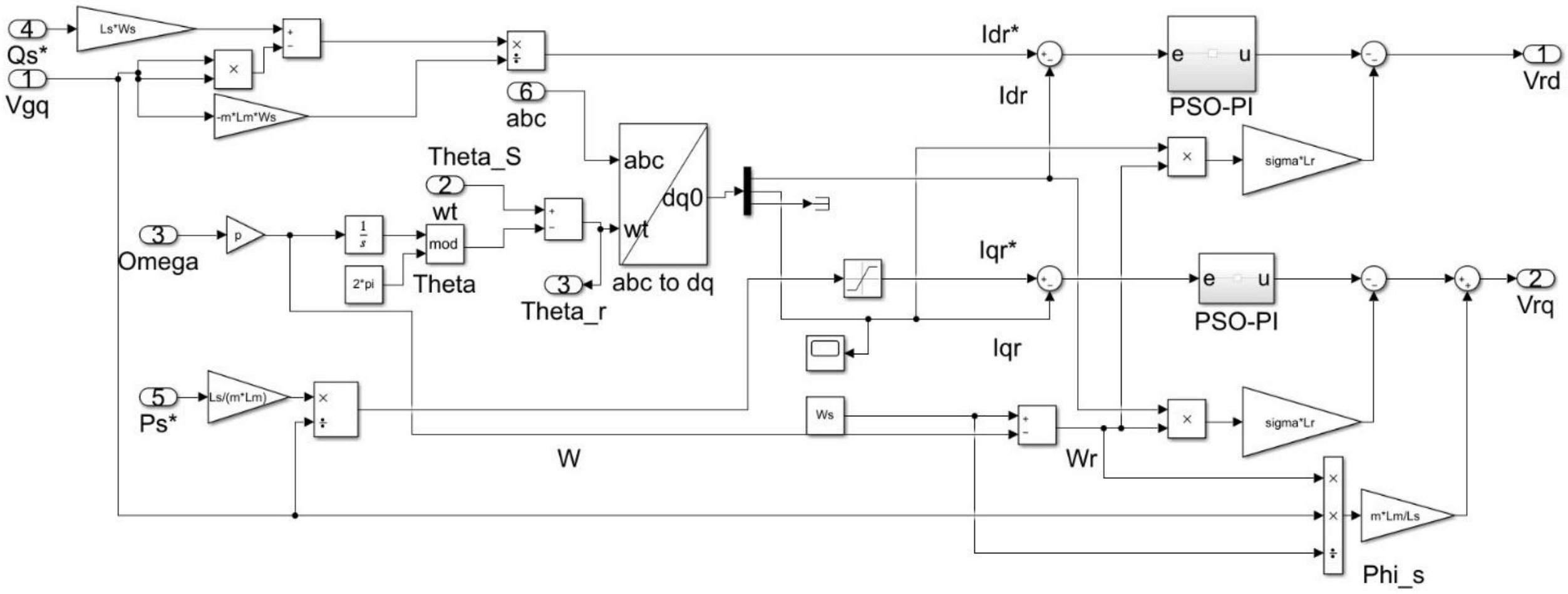

**Fig. 16: Simulink model of the DFIG control.**

The output voltages are then used to generate the PWM signal which controls the converter.

### 2.4.6. Filter and grid side converter (GSC)

As shown in Fig. 17, the GSC is linked to the grid through an RL filter, it allows the control of the DC bus voltage and the reactive power injected into the grid. Unity power factor is ensured by setting the reactive power to zero (Q=0).

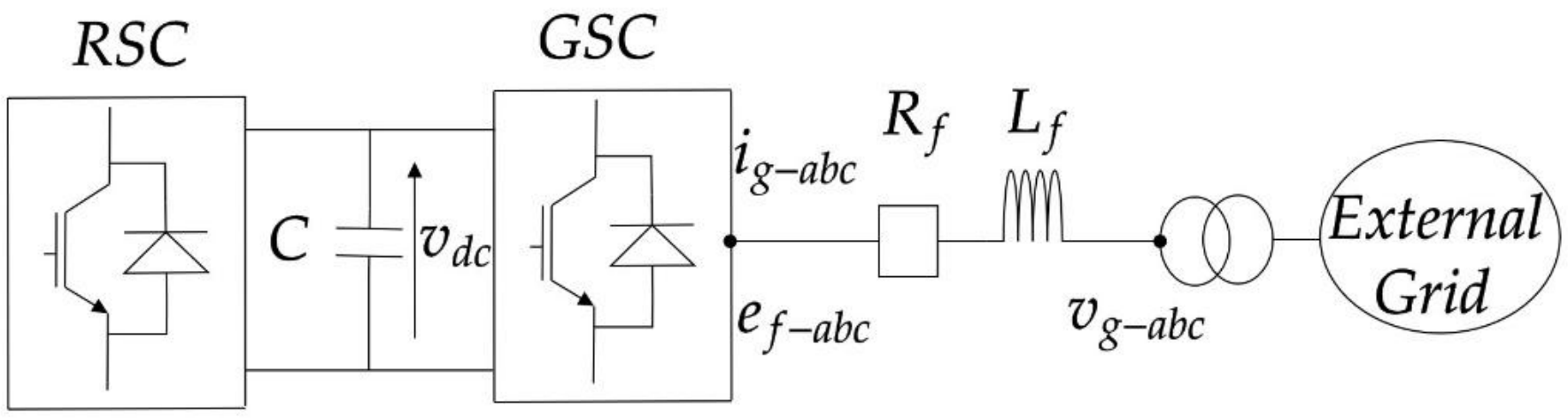

**Fig. 17: Connection between the back-to-back converter and the grid.**

The three-phase equations of the grid-side circuit are listed in equations (60-62):

$$v_{ga} = R_f . i_{ga} + L_f \frac{di_{ga}}{dt} + e_{fa} \quad (60)$$
$$v_{gb} = R_f . i_{gb} + L_f \frac{di_{gb}}{dt} + e_{fb} \quad (61)$$
$$v_{gc} = R_f . i_{gc} + L_f \frac{di_{gc}}{dt} + e_{fc} \quad (62)$$

Transforming the equations to the synchronous rotating frame, the d-q representation can be derived as:

$$v_{fd} = -R_f . i_{fd} - L_f \frac{di_{fd}}{dt} + e_{fd} \quad (63)$$
$$v_{fq} = -R_f . i_{fq} - L_f \frac{di_{fq}}{dt} - e_{fq} \quad (64)$$

Where:

$$e_{fd} = \omega_s . L_f . i_{fq} \quad (65)$$
$$e_{fq} = -\omega_s . L_f . i_{fd} + v_{gq} \quad (66)$$

The GSC allows the control of the RL filter currents $i_{fd}$ and $i_{fq}$ independently, the latter is done by PSO-PI controller based on the reference value $i_{fd}*$ and $i_{fq}*$. Knowing that the active and reactive powers from the GSC are given by:

$$P_f = V_{fd} . I_{fd} + V_{fq} . I_{fq} \quad (67)$$
$$Q_f = V_{fq} . I_{fd} + V_{fd} . I_{fq} \quad (68)$$

By neglecting the losses due to the RL filter resistance while considering that the voltage vector is aligned with the q-axis of the synchronous rotating frame, the d-q voltage components are obtained as:

$$v_{sd} = 0, v_{sq} = \widehat{v_g}.$$

The active and reactive powers and the reference active and reactive powers are given as:

$$P_f = V_{sq} . I_{fq} \quad (69)$$
$$Q_f = V_{sq} . I_{fd} \quad (70)$$
$$P_f^* = I_{fq}^* . V_{sq} \quad (71)$$
$$Q_f^* = I_{fd}^* . V_{sq} \quad (72)$$

The power in the DC bus is divided into:

$$P_{RSC} = V_{dc} . I_{RSC} \quad (73)$$
$$P_C = V_{dc} . I_c \quad (74)$$
$$P_{GSC} = V_{dc} . I_{GSC} \quad (75)$$
$$P_{RSC} = P_c + Joule \quad (76)$$

By neglecting power losses due to the joule effect in the capacitor, the converters and the RL fil become and $P_f$ becomes equal, and hence controlling $P_f$ means controlling $P_c$ which also translates into controlling the DC bus voltage ($V_{dc}$). As seen in the GSC Simulink model (Fig. 18), the dc-link voltage is controlled by the rotor q current component. Table 4 (see Appendix) lists the different parameters of the wind system.

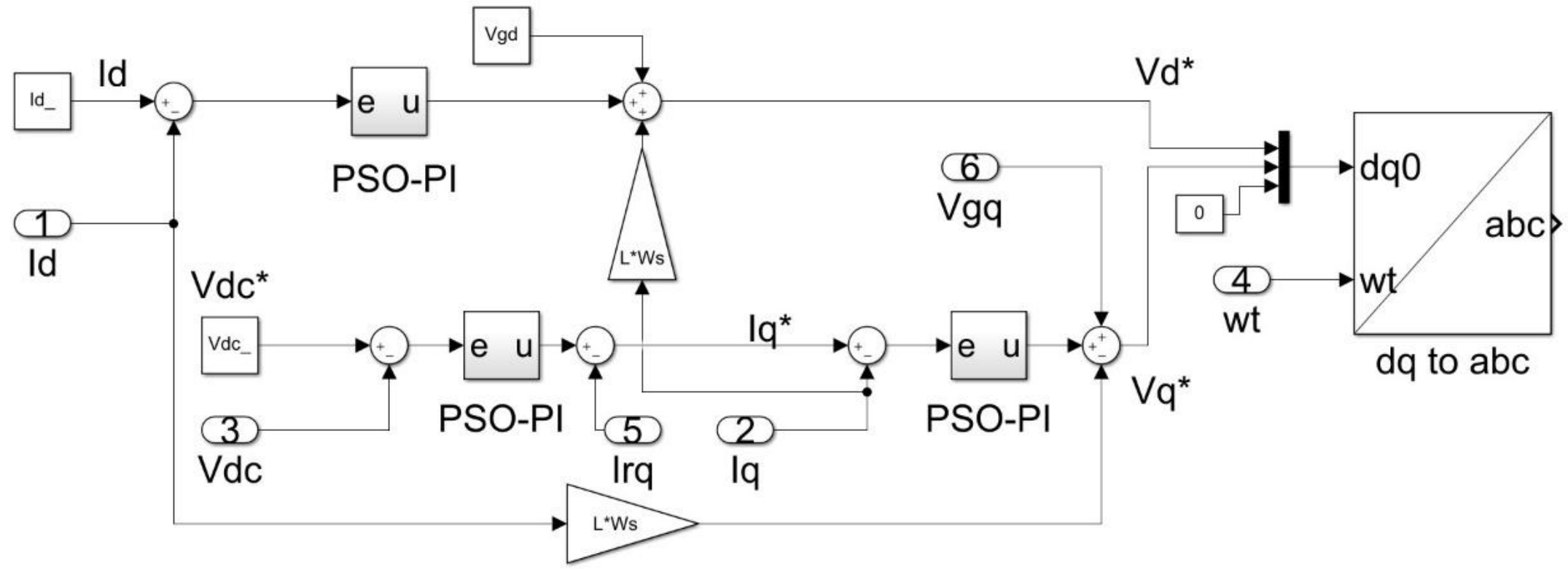

**Fig. 18: Simulink model of the GSC control.**

## 2.5. PSO-PI Control

PI controllers are widely used in the industry thanks to their simple design and good performance. The controller uses the error signal to generate a control command which is governed by the following equation:

$$u(t) = k_p e(t) + k_i \int e(t) d\tau \qquad (77)$$

where $k_p$ and $k_i$ are the proportional and integral gains and e(t) is the error. The discrete form of the PI has been used in this work where the control law is obtained using the trapezoid discretization method with $T_s$ as sampling time. It is often hard to properly tune the gains of the controller, therefore, the proportional and integral actions are optimized by an improved PSO. This meta-heuristic approach is a well-established technique. The algorithm relies on the principle of swarm social behavior, where the swarm contains several particles that represent potential solutions. Each particle is characterized by a position $p_i$ and a velocity $v_i$, the objective is to find the velocity that moves particles towards the optimal position. Further details can be found in (Kebbati et al. 2021), the algorithm is governed by the velocity and position update equations:

$$\begin{cases} v_i(k+1) = \omega v_i(k)) + c_1 r_1 (Pb_i(k) - x_i(k))+) + c_2 r_2 (Gb_i(k) - x_i(k)) \\ x_{i(k+1)} = x_{i(k)} + v_i(K+1) \end{cases} \qquad (78)$$

where $\omega$ is known as inertia weight, $c_1$ and $c_2$ are known as cognitive and social accelerations coefficients respectively. $r_{1,2} \in [0,1]$ are random constants, and $Gb$ and $Pb$ represent the global best position of the whole swarm and the local best position in the current swarm generation. In the classic PSO, the inertia weight, social and cognitive accelerations are prefixed constants. However, in the improved version of this work they are dynamic and change with iterations according to equations (3.76,3.77), this enhances the PSO search performance:

$$\omega = \omega_{min} + \frac{\exp(\omega_{max} - \lambda_1(\omega_{max} + \omega_{min})\frac{g}{G})}{\lambda_2} \qquad (79)$$

$$\begin{cases} c_1(k+1) = c_1(k) + \alpha \\ c2(k+1) = c_2(k) + \beta \\ \alpha = -2\beta = 0.085 \quad for\ \frac{g}{G} \leq 30\% \\ \alpha = \frac{-\beta}{2} = 0.045 \quad for\ 30\% \leq \frac{g}{G} \leq 60\% \\ \alpha = -\frac{\beta}{2} = -0.025 \quad for\ 60\% \leq \frac{g}{G} \leq 85\% \\ \alpha = -\beta = -0.0025 \quad for\ \frac{g}{G} \geq 85\% \end{cases} \tag{80}$$

The terms $g$ and $G$ represent the actual and the last generations respectively, $\lambda_{1,2}$ are constant parameters adjusted to ensure an exponential decrease from $\omega_{max}$ to $\omega_{min}$ which are maximum and minimum inertia weights. The advantage of this improved version compared to other versions is that it increases the overall search capabilities of the PSO algorithm, as the exponential decrease of $\omega$ accelerates the convergence towards the global best solution. Increasing $c_1$ pulls the particles towards $Pb$ and enhances the exploration phase, while increasing $c_2$ speeds the convergence towards $Gb$ which enhances the exploitation phase and vice versa. The proposed PSO is used to optimize the PI actions $\{k_p, k_i\}$ applied to the complete model. The approach is illustrated in Fig. 19 where the mean squared error (MSE) is used as the fitness function for the PSO algorithm. The model is simulated with the gains obtained from the PSO algorithm, which are optimized at every iteration until the optimal solution is found. (Kebbati et al. 2021) applied the PSO algorithm to optimize a model predictive controller for path tracking applications in autonomous vehicles. The same algorithm was implemented in (Kebbati et al. 2022) to tune a PID controller for autonomous driving applications.

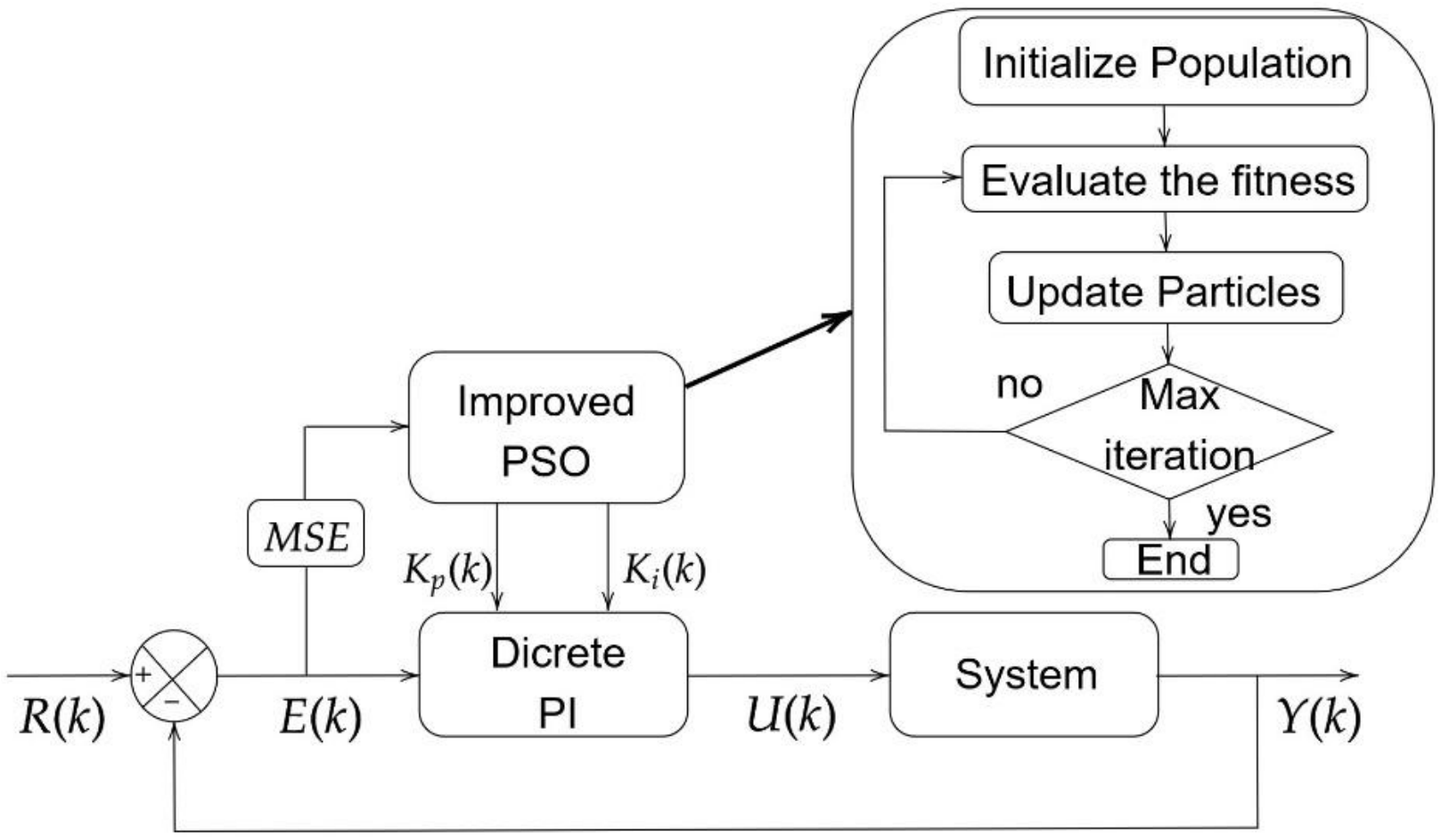

**Fig. 19: PSO-PI control approach**

The optimization with improved PSO resulted in the gains $\{k_{p1} = 5.6749, k_{i1} = 11.6077\}$ for the inverter controller, and the gains $\{k_{p2} = 3.4074, k_{i2} = 9.4171\}$ for the DFIG controller, while the optimized gains for the grid side converter were found as $\{k_{p3} = 7.6999, k_{i3} = 5.019\}$.

## 2.6. Economic and Environmental Analysis

In this study, the hybrid system is of large scale and injects power directly into the grid. Such systems have very high initial investments but pay off in the long run. The size of the system is important in determining the feasibility and cost-eff of the whole project. In general, the cost of any system is composed of capital costs, replacement, operating and maintenance costs. The capital cost must include the price for all the components of the system as well as transportation, packaging and installation of these components, etc. The capital cost is given by the following formula:

$$Capital\ cost = (1 + K).N.C \tag{81}$$

Where C and K are the cost of one component (PV panel, wind turbine, etc.) and the cost related to engineering, logistics, installation, etc. N represents the number of components (100 PV panels, 30 converters, etc.). On the other hand, the replacement and operating and maintenance costs are known as variable costs. The operating and maintenance costs include the expenses of insurance, salary of operators, taxes, recurring costs and maintenance. They are expressed as a percentage of the initial capital cost. The replacement costs are the expenses needed to change components as these have different lifetimes and some of them requires to be replaced several times during the lifespan of the system. When conducting economic analyses of energy systems, the main factor for judging if the system is feasible or not is the levelized cost of electricity (LCOE). The latter represents the cost of generated energy units during the lifespan of the system. The LCOE is calculated by summing the capital costs and operating and maintenance costs and dividing by the total energy generated from the system. The net present value (NPV) approach is used for the calculation of the LCOE. In this method, the total expenses and the incomes of the system are determined based on discounting from a starting date, the LCOE is given as:

$$LCOE = \frac{I_0 + \sum_{t=1}^{n} \frac{A_t}{(1+i)^t}}{\sum_{t=1}^{n} \frac{M_{t,el}}{(1+i)^t}} \tag{82}$$

Where $I_0$, $A_t$ and $M_{t,el}$ are the respective investment expenditure, annual total cost and produced quantity of electricity. $i$, $n$ and $t$ represent the interest rate, the economic operational lifetime and the year respectively. The NPV measures the net profit and actual value of the project, it includes all benefits and cost streams occurring at different points in time and converts them into present value equivalents. The NPV is calculated over a period of time by subtracting cash outflows from the values of cash inflows. The overall worth of the project is determined by the present value equivalents, the formula for NPV is given as:

$$NPV = \sum_{t=1}^{n} \frac{R_t - C_t}{(1+i)^t} - I_0 \tag{83}$$

Where $R_t$ and $C_t$ represent the revenue and the costs in year t, $i$ and $I_0$ are the discount rate and the initial investments respectively. The payback time represents the time required for the system to amortize the initial project investment by its cash inflows. Shorter payback times ensure a faster recovery of the investment funds. Generally, the distribution of cash flow is not considered beyond the payback period. The payback time can be calculated using the following equation:

$$Payback\ period = \frac{Initial\ Investment}{cash\ inflow\ per\ period} \tag{84}$$

### 2.6.1. System model in HOMER

Homer is used to perform an economic and environmental analysis of the system, which is supposed to inject all its power into the grid. Thus, all produced electricity will be sold to the utility. The different system components are shown in Fig. 20 and listed with their respective costs in Table 5. The cost of engineering, logistics and installation has been estimated at 25% of the component's cost. The interest rate is set to 3.75% while the inflation rate is evaluated at 2.42% (Ministry of finance, 2021). Moreover, the electricity price is fixed at 0.04 $/kWh. Electricity producers using solar and wind energies benefit from a premium of 300% of the market price of electricity, this is ensured for the full project lifetime according to law “No.° 02-01 du 22 Dhou El Kaada 1422 (Articles 88 ff)” (Meyer-Renschhausen. M. 2013).

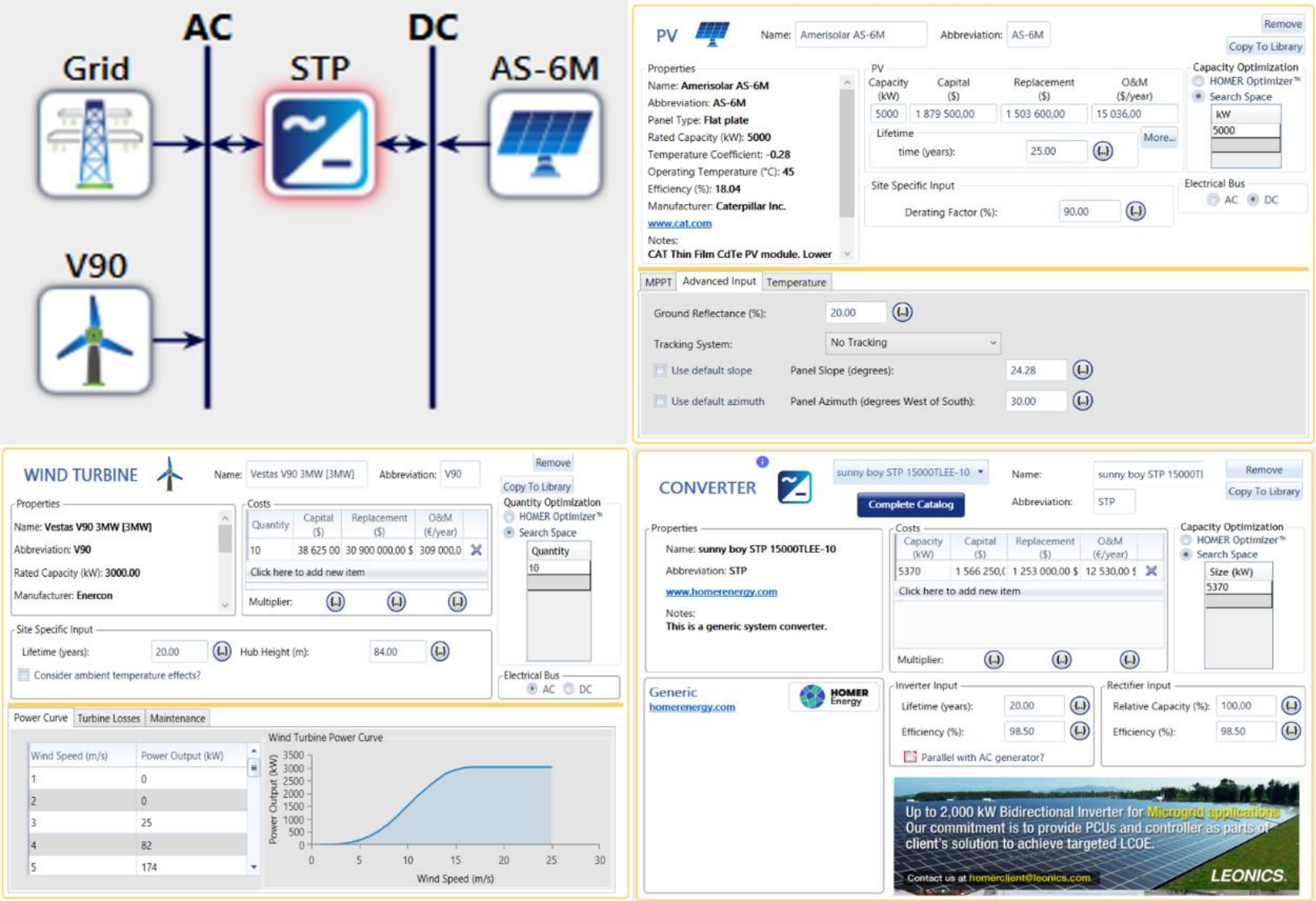

**Fig. 20: Homer model of the hybrid system**

For the environmental analysis, the amount of sulfur dioxide and nitrogen oxides emitted by a combined cycle power plant using natural gas have been chosen as 2.74 kg/MWh and 1.34 kg/MWh, while the amount of CO2 is estimated at 400 kg/MWh (Benhadji Serradj et al. 2021).

**Table 5: Summary of the system costs**

| Component | Capital cost ($) | Replacement cost ($) | OMC (%) |
|---|---|---|---|
| PV panel (AS-6M) | 132 | 105 | 1 |
| Wind turbine (Vestas V90) | 3862500 | 3090000 | 1 |
| DC/AC converter (Sunny boy STP 15000) | 4375 | 3500 | 1 |

# 3. Results and discussion

The system is simulated for two specific months using hourly data that was generated by METEONORM. The result is hourly power output; the latter is compared to hourly power production of the central power plant of Adrar. The inputs to the PV system are global horizontal irradiation, direct normal irradiation, diffuse horizontal irradiation and the number of hours. The outputs are active and reactive powers. On the other hand, the wind system receives hourly wind speeds, and outputs active and reactive powers as well. The performance of the control is presented along with the energy output of the PV plant and the wind farm. An analysis of the amount of power injected into the grid and the contribution of the hybrid system is also conducted. The performance output of the PV system is affected mainly by the cell's temperature and the irradiation received by the PV generator, the temperature lowers the efficiency. The trend of cell temperature levels in January shows slight variations. Peak levels vary each day with a maximum of 52 °C at the beginning and the end of the month, and a minimum of 41 °C is noticed around mid-January. Fig. 21 shows hourly variations of cell temperature in January versus August.

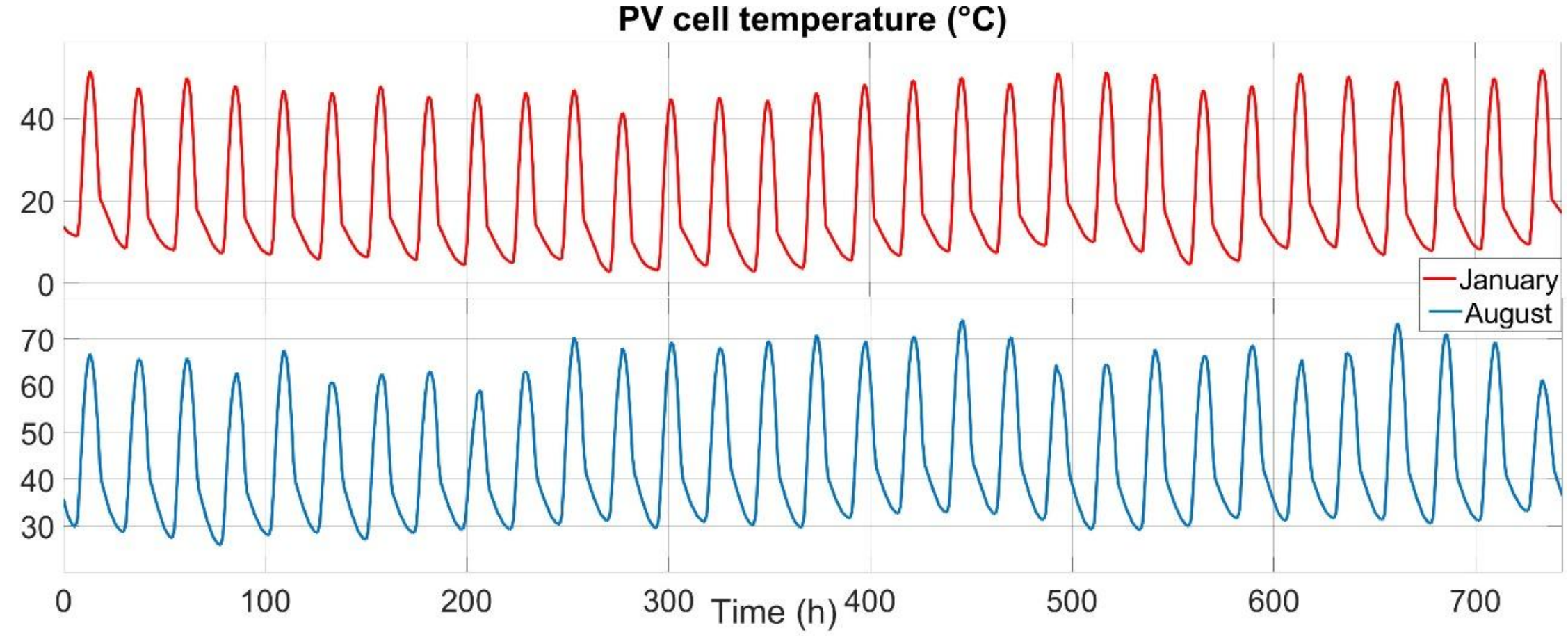

**Fig. 21: PV cell temperature variations (January vs August).**

Although ambient temperatures in January are significantly low, the irradiation received by the module is unpredictably high (see Fig. 22). Minimum radiation levels of around 850W/m² are observed in the first and last days of January, this can be attributed to the solar position throughout the year as the distance between the earth and the sun varies and leads to a variation of extraterrestrial radiation flux. Irradiation levels in August are considerably low compared to January, low irradiation coupled with high temperatures reduces dramatically the PV output.

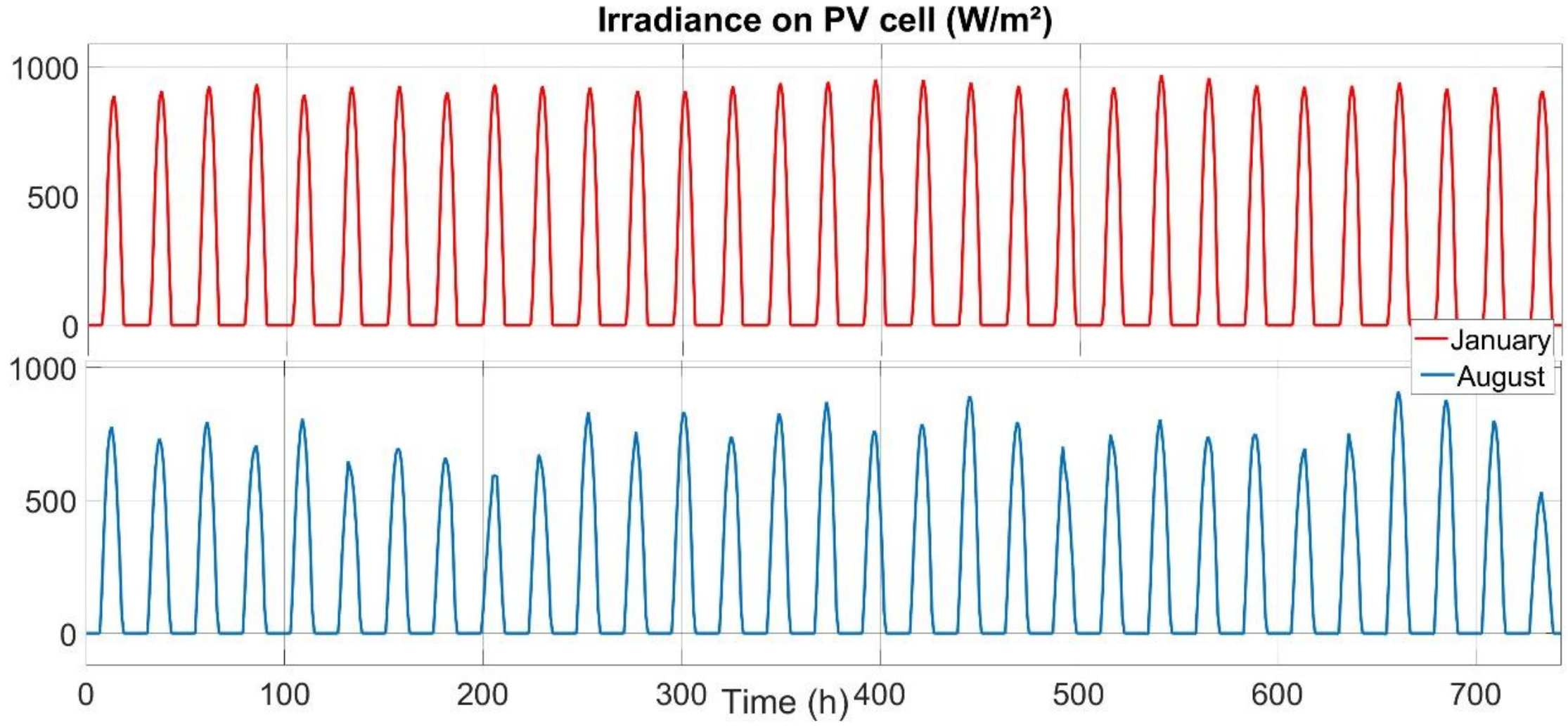

**Fig. 22: Irradiance received by the PV modules (January vs August).**

### 3.1. MPPT control performance for PV system

The perturb and observe algorithm keeps the voltage output of the PV module at its maximum value, the accuracy is determined by the size of the perturbation (ΔV) and thus the output keeps oscillating around the maximum value as demonstrated by Fig. 23. It is worth mentioning that the MPPT tracker is fast, and the maximum power point is quickly reached, hence the control works perfectly.

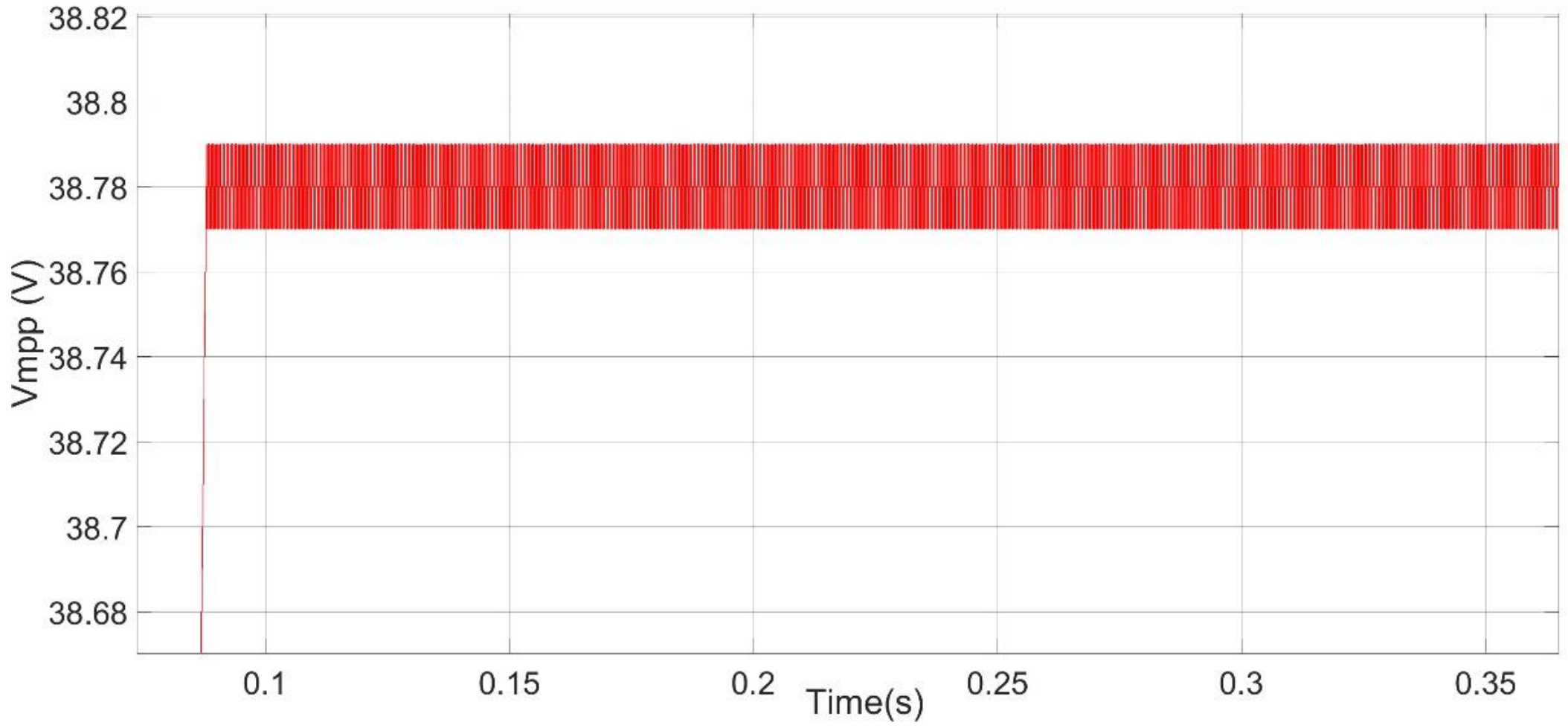

**Fig. 23: Performance of the MPPT tracker.**

Fig. 24 shows how the power output of the PV module follows the maximum power point under standard test conditions (STC) (G=1000 W/m² and $T_c$=25°C).

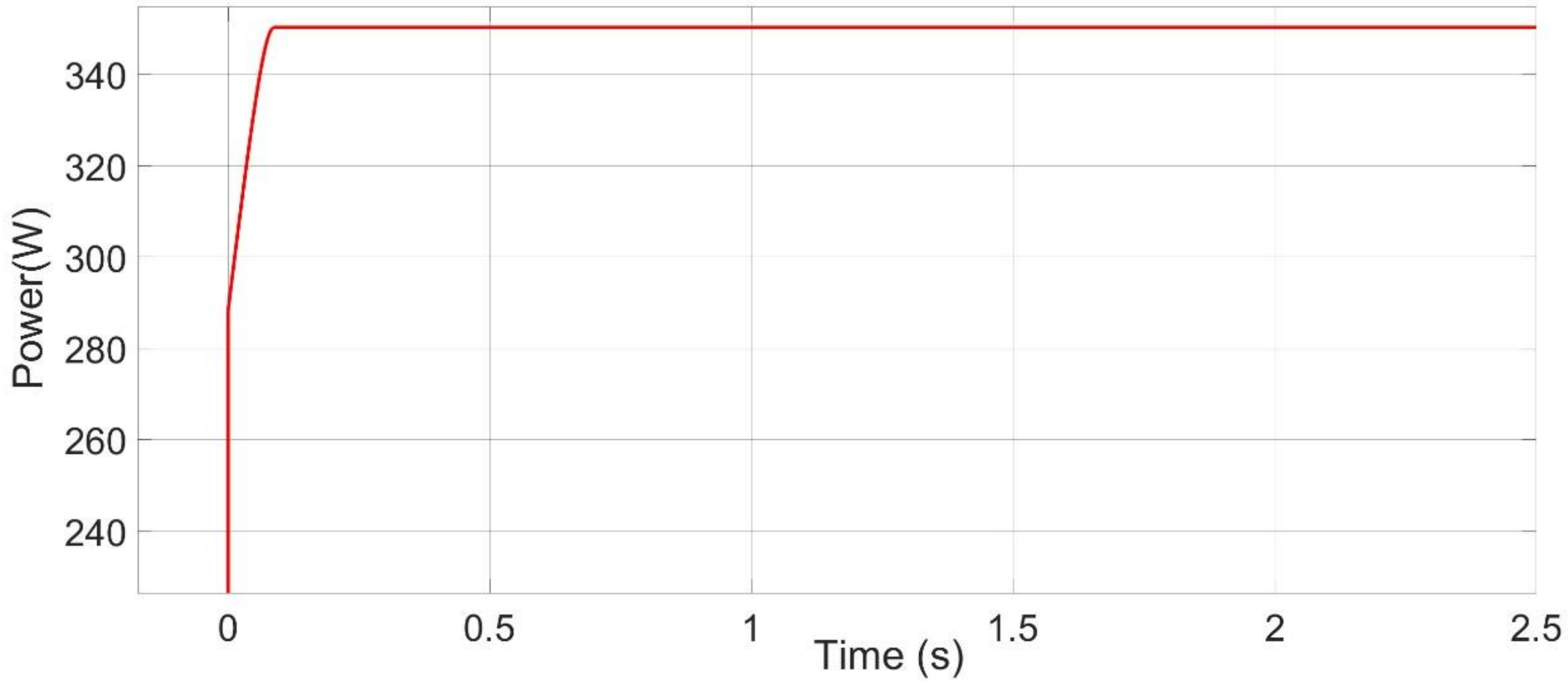

**Figure 24: PV panel power output at STC.**

### 3.2. DC-DC converter control

The boost converter control maintains $V_{dc}$ constant all the time as observed in Fig. 25, and the latter will be received by the inverter to be converted into AC form.

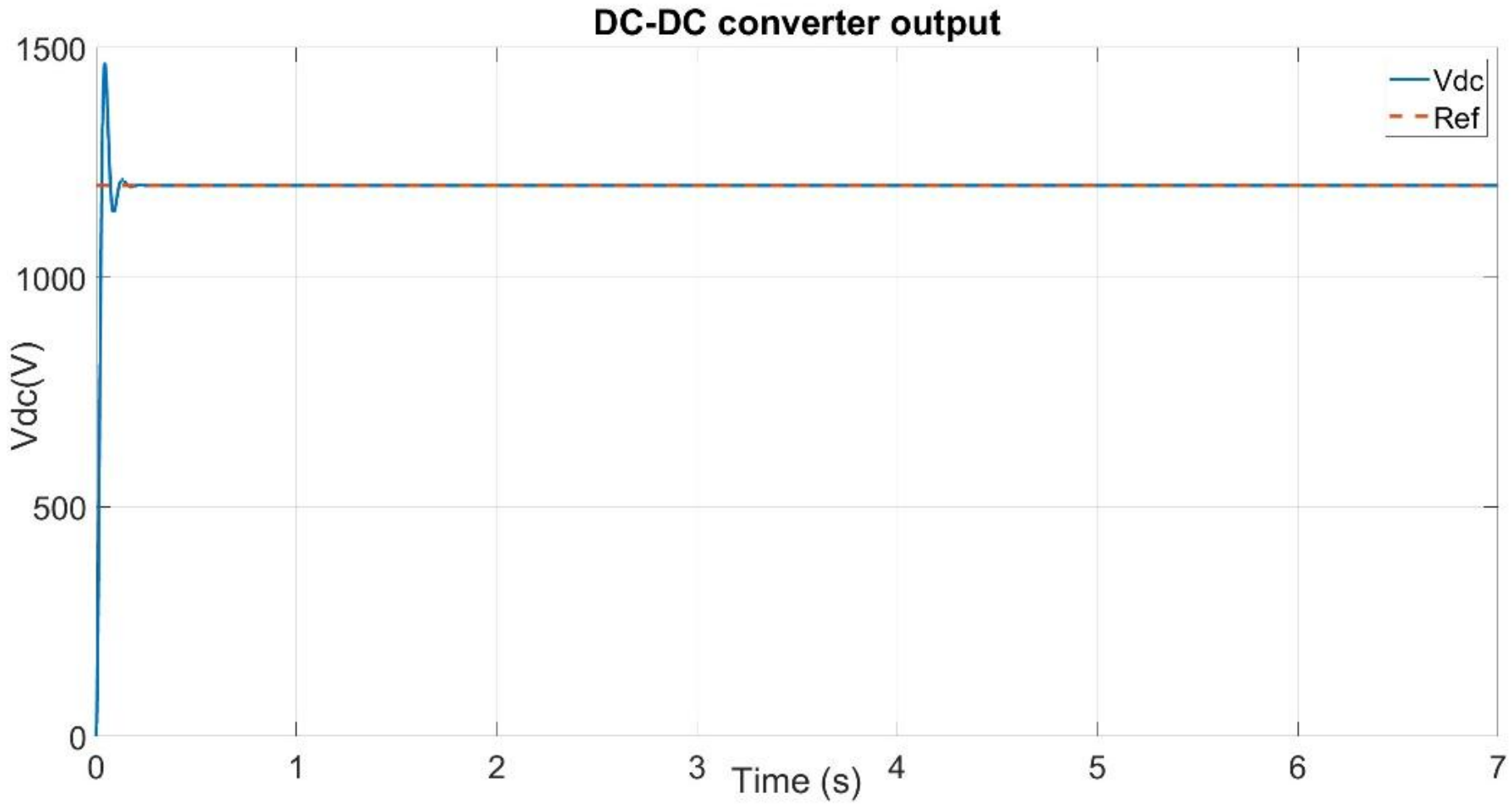

**Fig. 25: DC-DC Converter output.**

## 3.3. Active and reactive power control for PV system

The performance of the PV system power control is shown in Fig. 26, both active ($P$) and reactive ($Q$) powers follow exactly the reference values. Negative values imply that the energy is injected into the grid while positive values mean that the energy is drawn from the grid. The reactive power is kept around zero to ensure unity power factor. The amount of power generated corresponds to one string under STC conditions, which is equivalent to 14 KW.

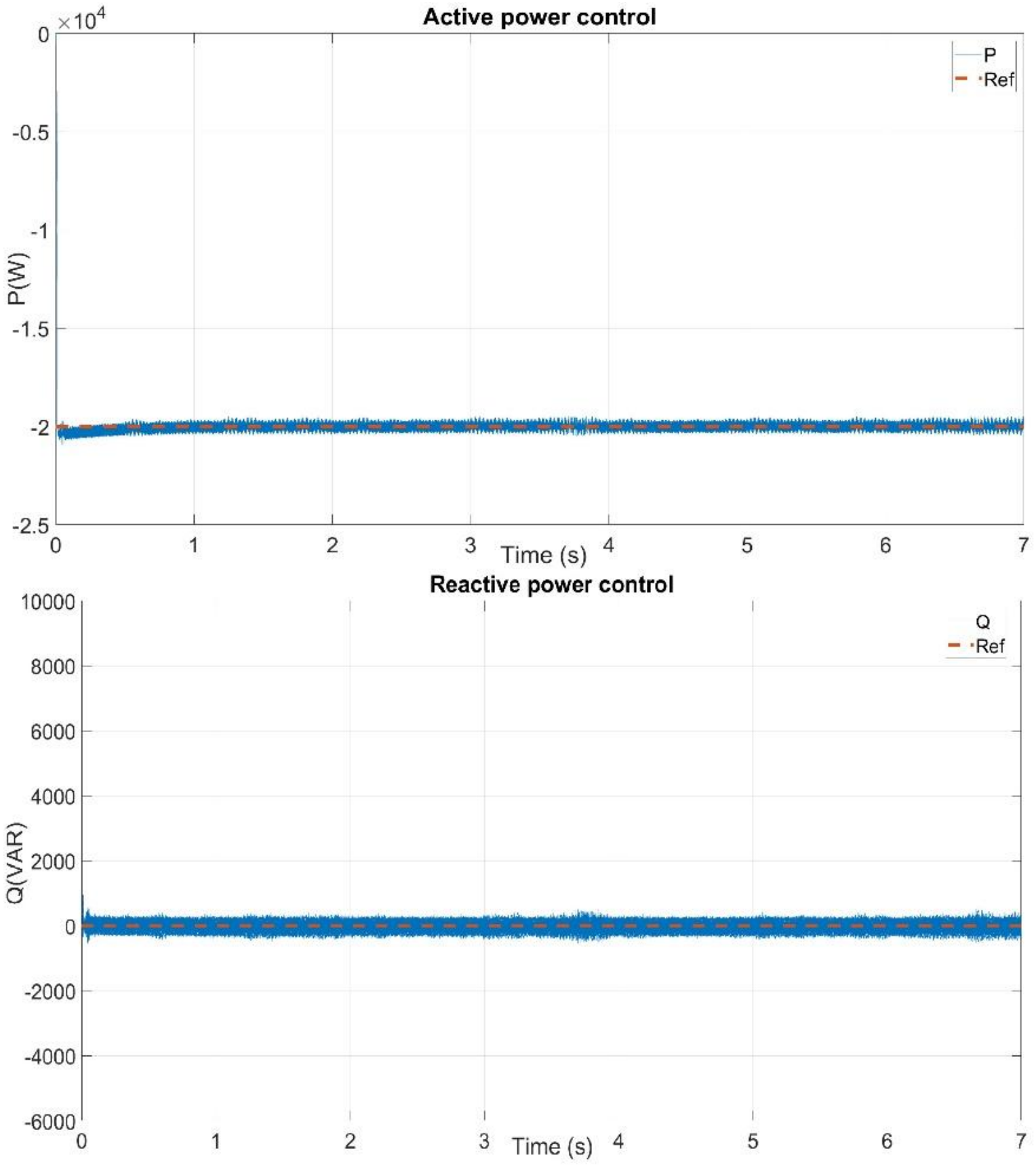

**Fig. 26: Power control for the PV system**

### 3.4. Output of the PV plant

The simulation results for January show that the PV generator produced from 4.35 MW to 4.73 MW depending on the day, irradiation and temperature levels. The power is produced from 7 am to 6 pm with peak production occurring between midday and 2 pm. During the month of August, the PV generator produces less power compared to January. This decrease in power output is justified by the high temperatures which reduce the efficiency of the PV generators. Peak production ranges from 2.45 MW to 4.15 MW. However, the production during the day lasts from 7 am to 7 pm delivering power for an extra hour compared to January, since summer days are longer than winter days.

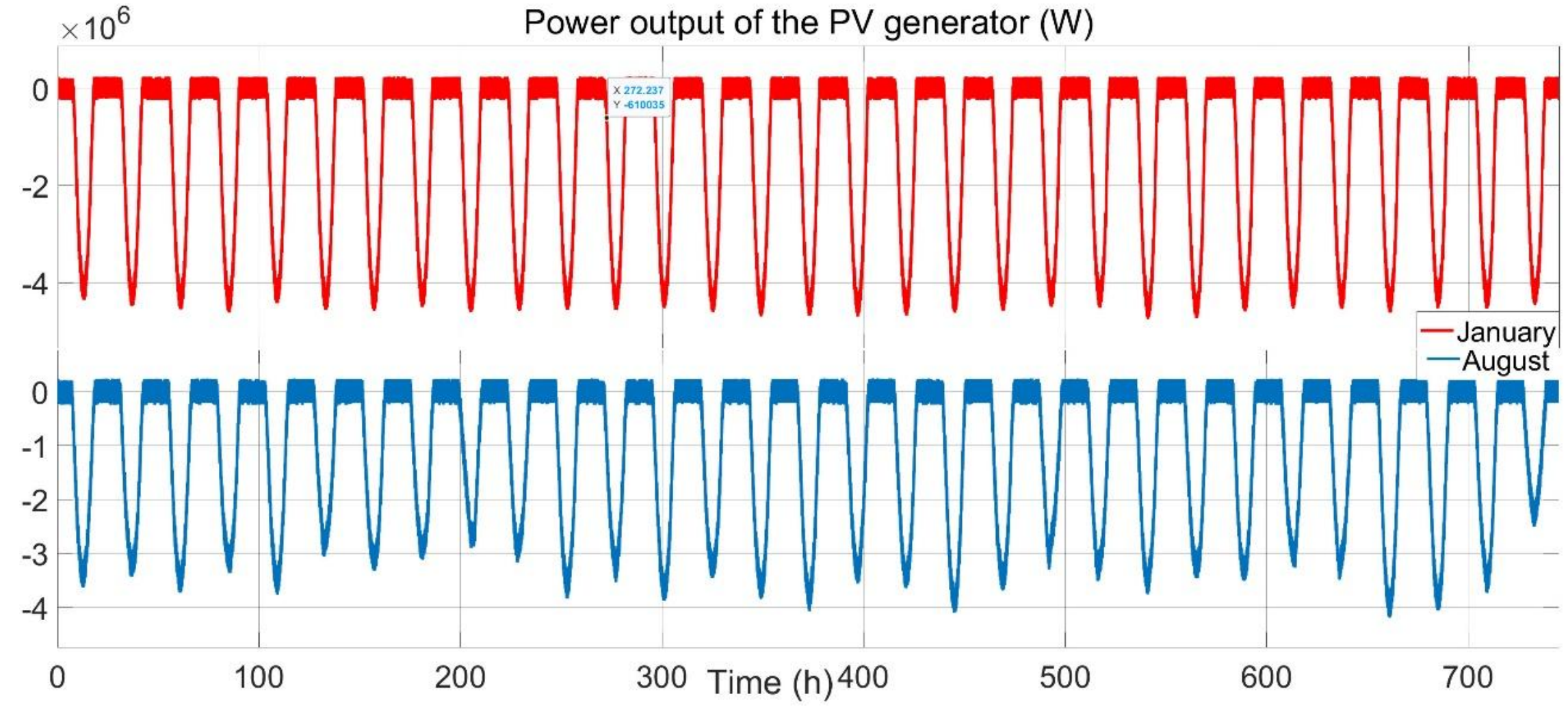

**Fig. 27: Power output of the PV plant (January vs August).**

The wind system has a rated power of 30 MW, which depends on the speed of the wind. The turbine's pitch angle is fixed at its optimal value and the hub height is fixed at 80 meters, where the corresponding wind speeds are shown in Fig. 28. In August, the wind speed profile shows more energy potential as high speeds occur more frequently. The average speed is around 8 m/s, meaning that the system should generate more energy in August than in January. Thus, the performance of the grid should be improved in this period of time that corresponds with high demand.

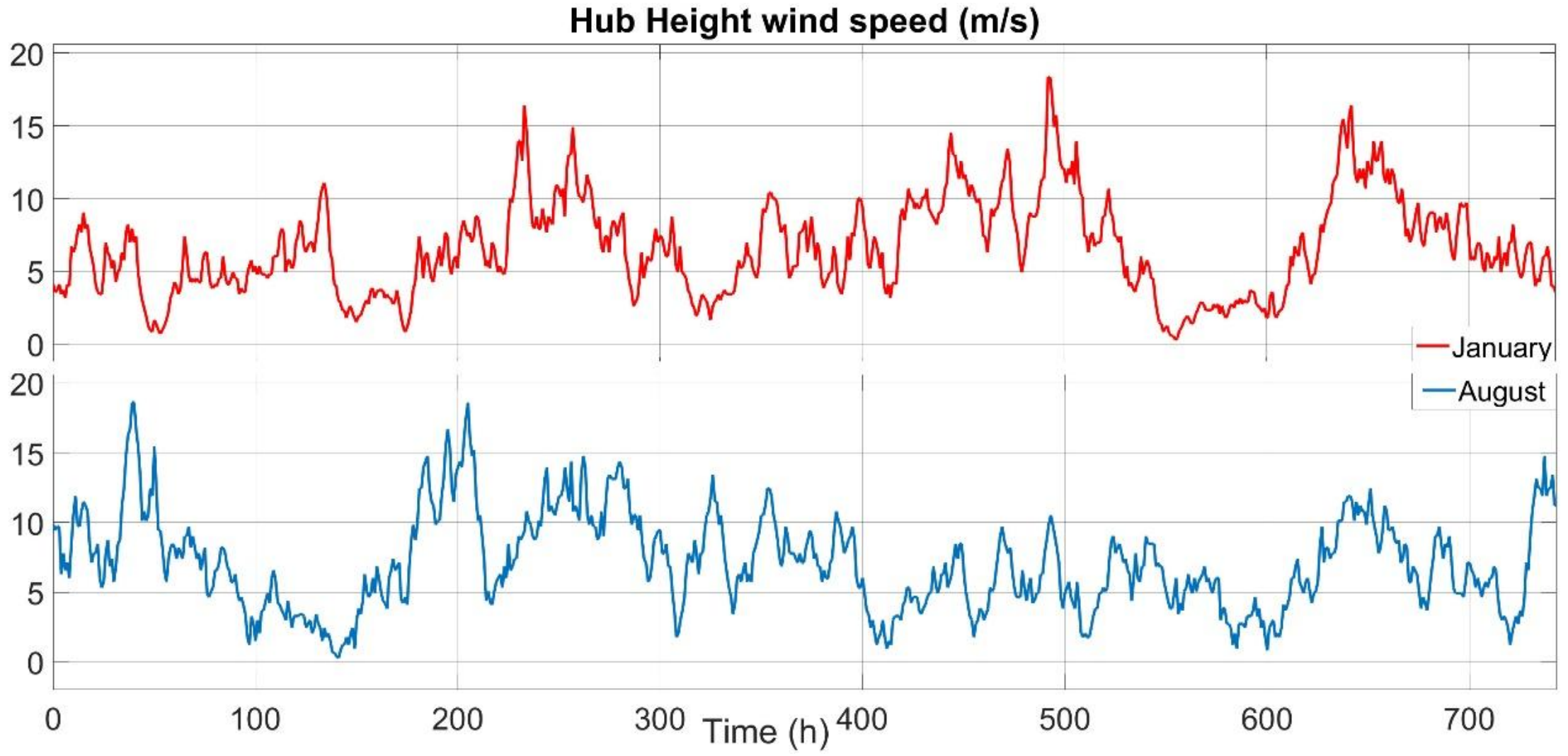

**Fig. 28: Wind speeds at hub height (January vs August).**

### 3.5. MPPT control performance for wind system

The MPPT implemented to control the rotation speed of the turbine was tested with a multitude of different reference values and the result is presented in Fig. 29. The PSO-PI performs very well ensuring maximum power extraction from the wind. As observed, the rotational speed follows exactly the reference values, and the robustness of the control is tested by changing the reference values every time. Fig. 30 shows the performance of the wind system power control. In the first graph of Fig. 30, the active power follows exactly the given reference values. The reference values in the system are obtained from the mechanical model as explained in the control section. In the second graph of Fig. 30, the reactive power is maintained at zero which ensure a unity power factor.

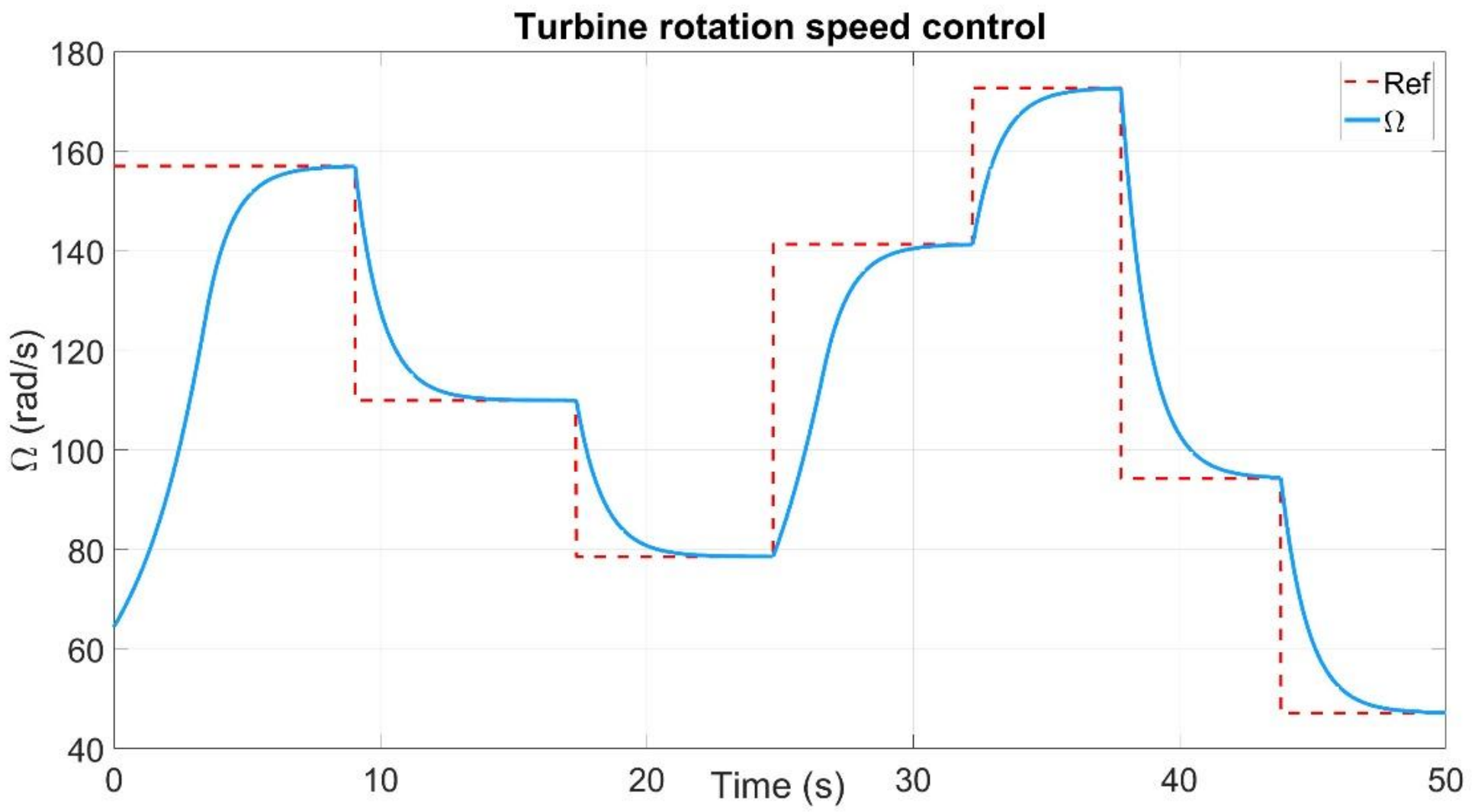

**Fig. 29: MPPT turbine rotation speed tracking.**

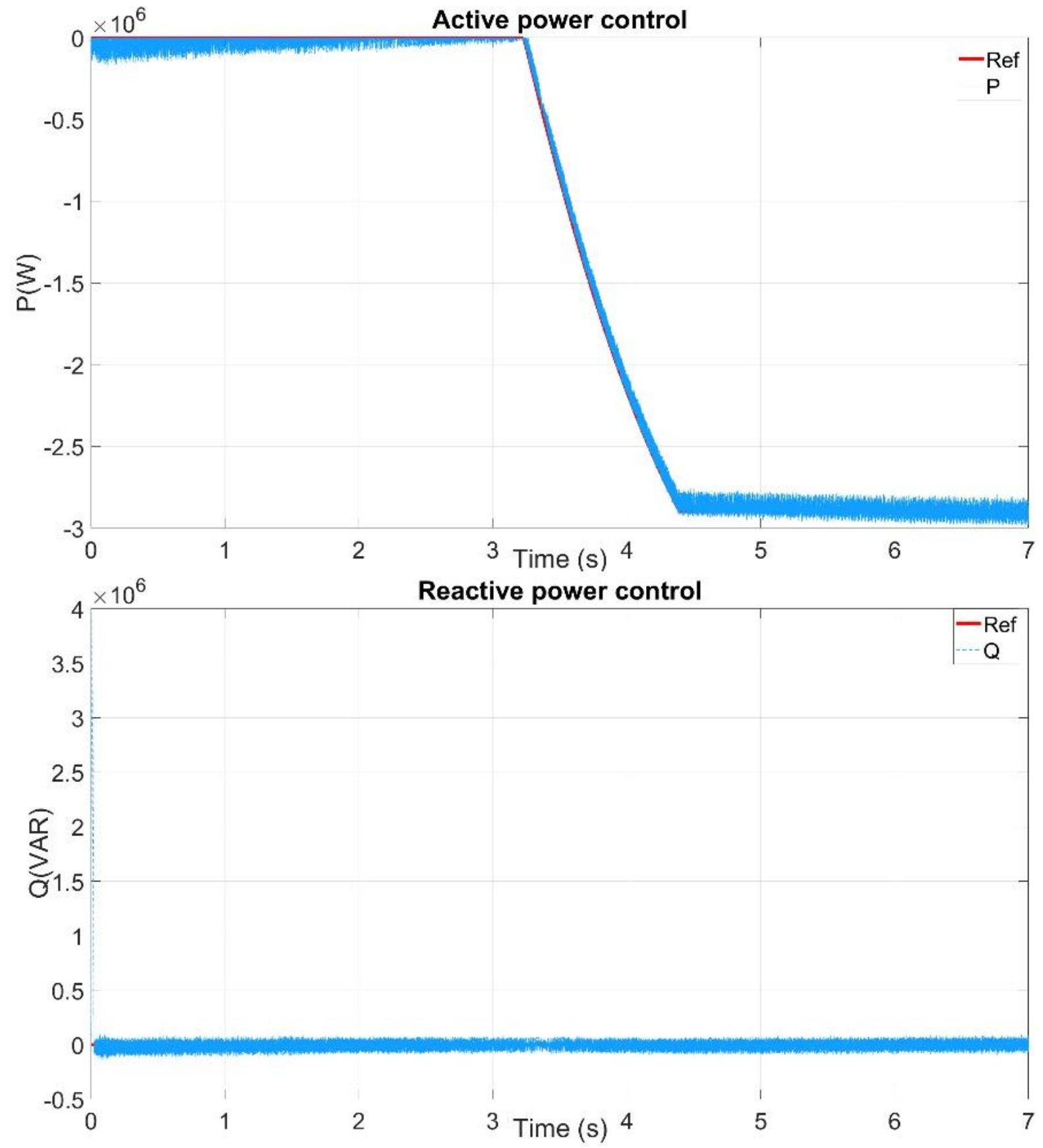

**Fig. 30: Power control for the wind system.**

### 3.6. Output of the wind farm

January is characterized by lower wind speeds; hence the power output is lower compared to August. Peak production is reached in three weeks out of four with an output ranging from 10 to 23 MW. During August, more power was produced and the rated output was often reached especially in the beginning and at the end of the month. An average of around 10 to 25 MW is maintained throughout the month. This amount of power will contribute significantly to grid performance as the cooling load increases during this month.

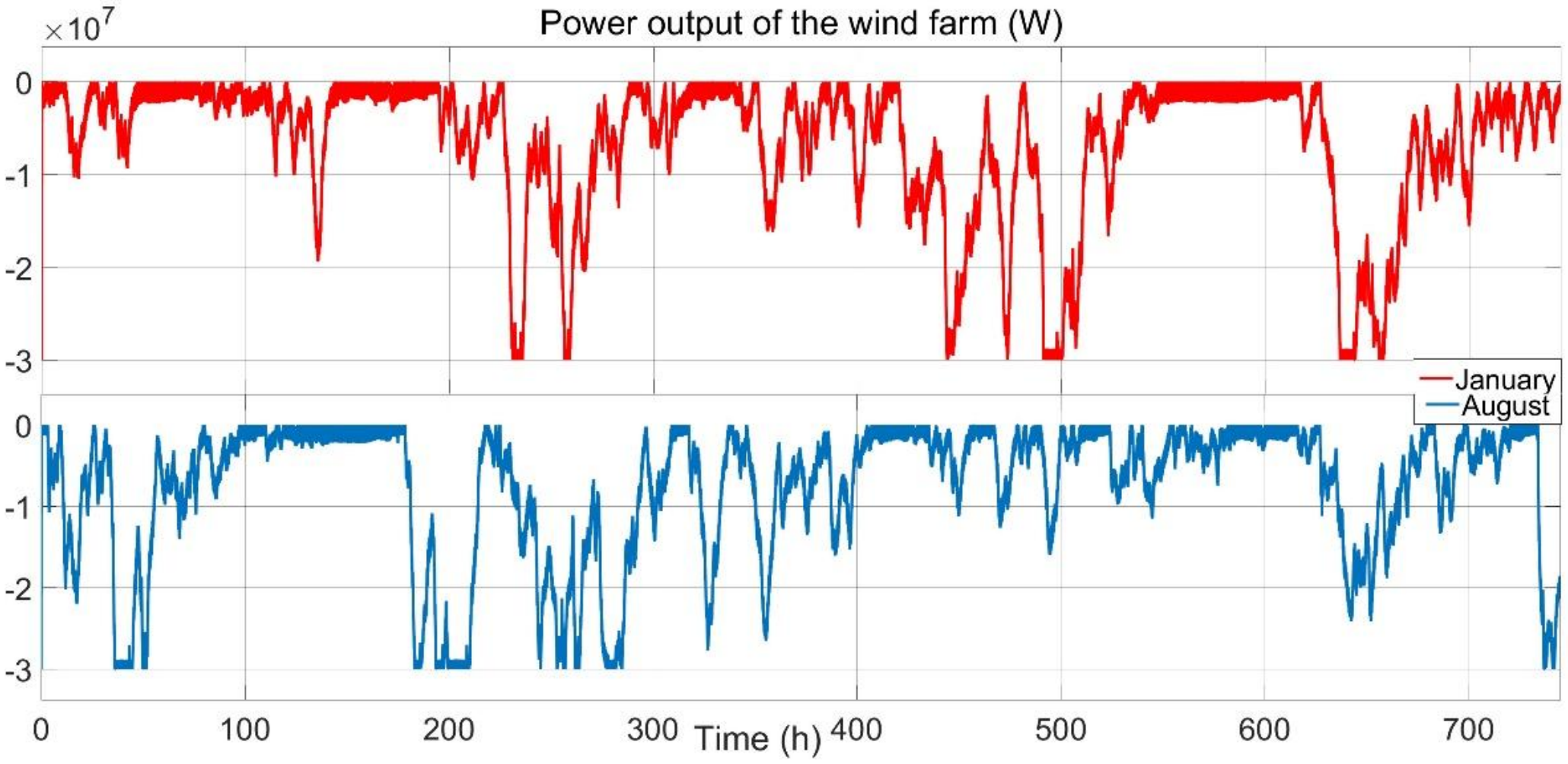

**Fig. 31: Power output of the wind farm (January vs August).**

## 3.7. Results of the hybrid system

The power produced by the hybrid system is compared to that produced by the central power plant of Adrar to determine the contribution of the system to the grid. More contribution is achieved in January varying from 10 to 34.5 MW while less contribution was achieved in August due to higher energy demand. The peak loads seemed to appear from 2 pm to 9 pm, while the lowest production occurs between 2 am and 10 pm. Despite the high power produced by the hybrid system, it is still low vis-a-vis the central power plant. However, on most days the hybrid system managed to cover up to 50% of what the plant produces.

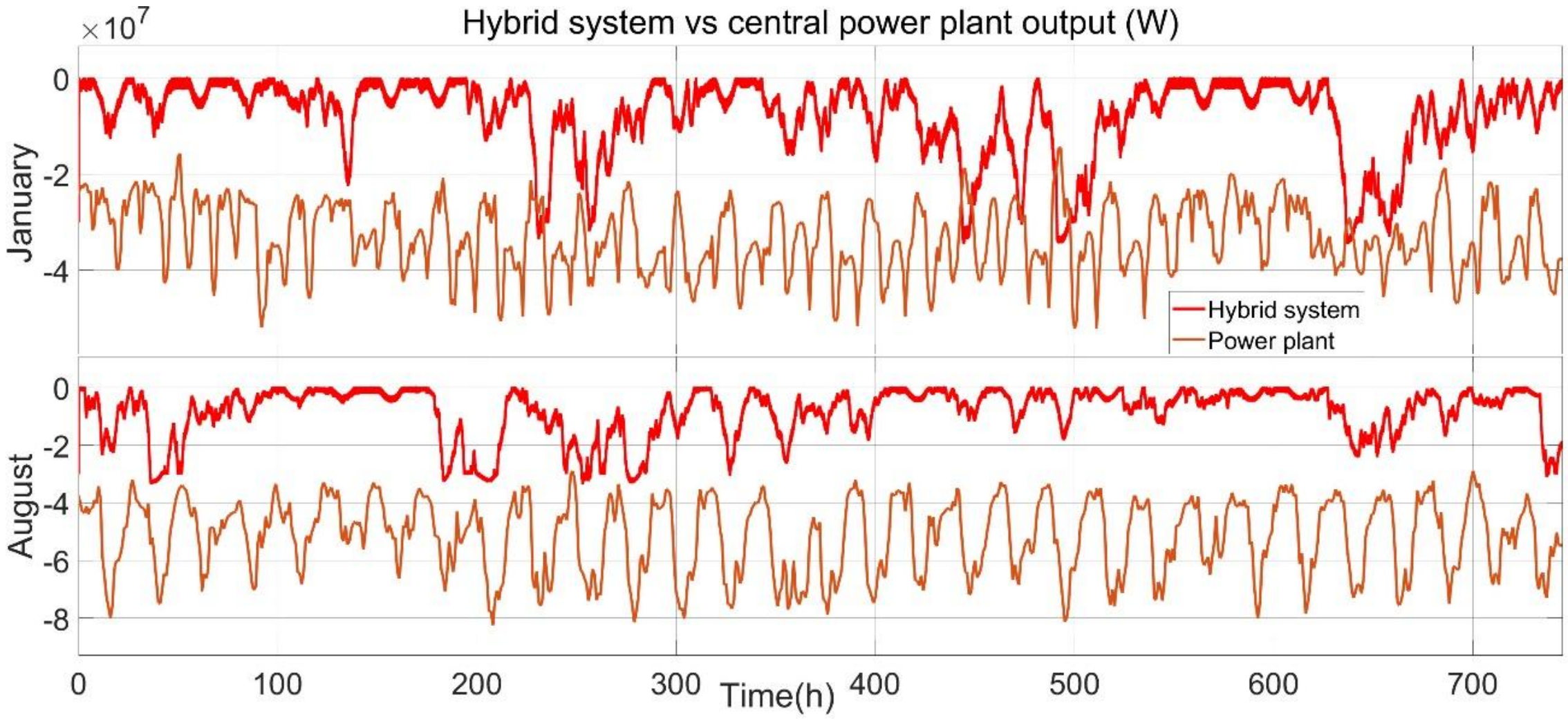

**Fig. 32: System power output vs central plant power production (January vs August).**

Fig. 33 shows the amount of power covered by the hybrid system. In January, the hybrid system covered from 50% to 220% of what the central power plant produced. While in August, the share covered by the system varies from 20% to 85%. Note that the central power plant supplies power to the whole province of Adrar, in reality the true load of Adrar city is much lower.

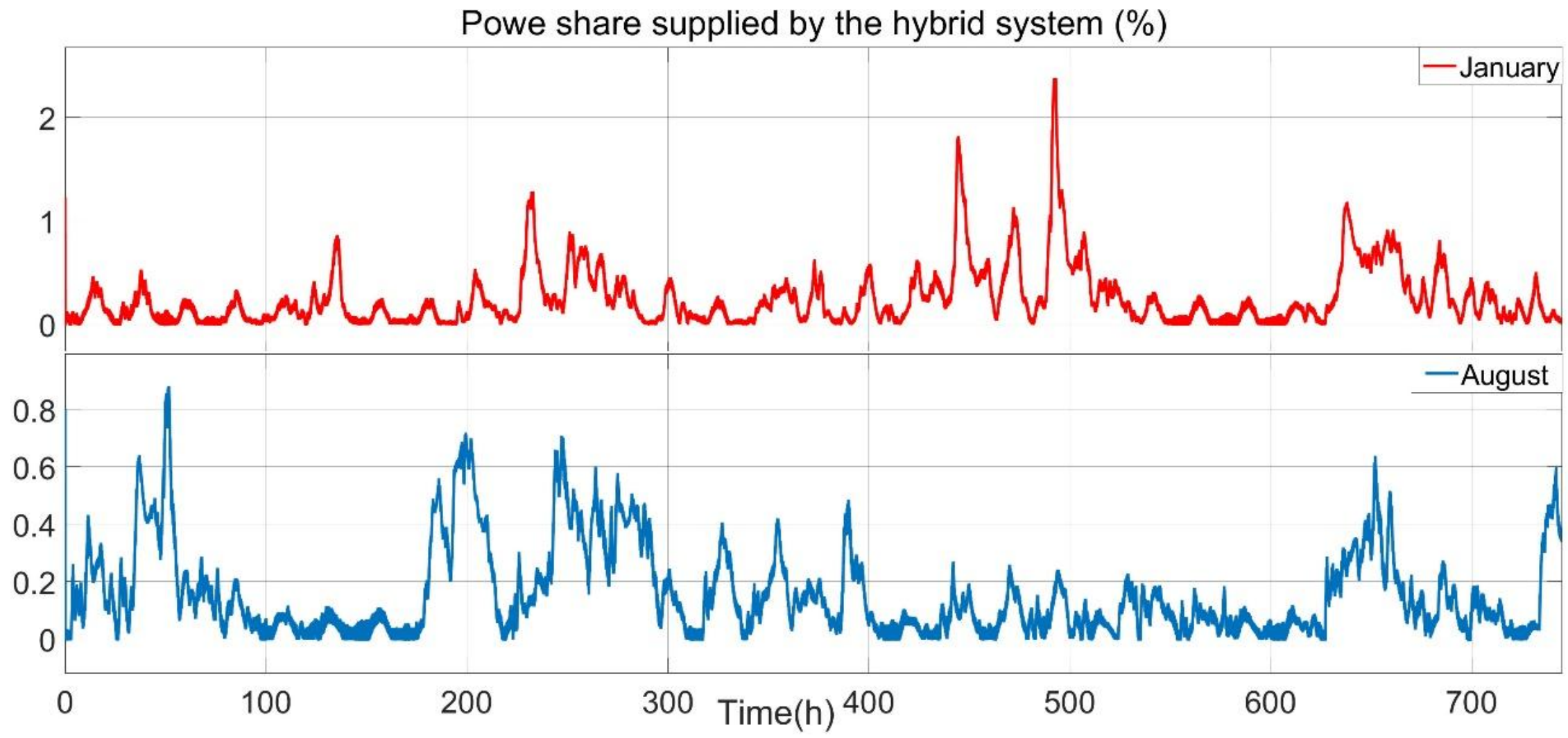

**Fig. 33: Share covered by the hybrid system (January vs August).**

### 3.8. Cost Analysis

HOMER results show from an economical point of view that the system composed of a 5 MW rated PV system and a 30 MW rated wind system has a total net present value (NPV) of $167 million, whereas the levelized cost of electricity (LCOE) sits at 0.084 $/kWh. The LCOE is very profitable compared to the selling price of 0.114 $/kWh guaranteed by the feed-in tariff for the duration of the project. As a result, the payback period is in just 7.3 years and by the end of the projects lifetime a cumulative net profit of $127 million is reached. Fig. 34 depicts the cash flow over the lifespan of the project, which shows how much money is either spent on or gained by the project.

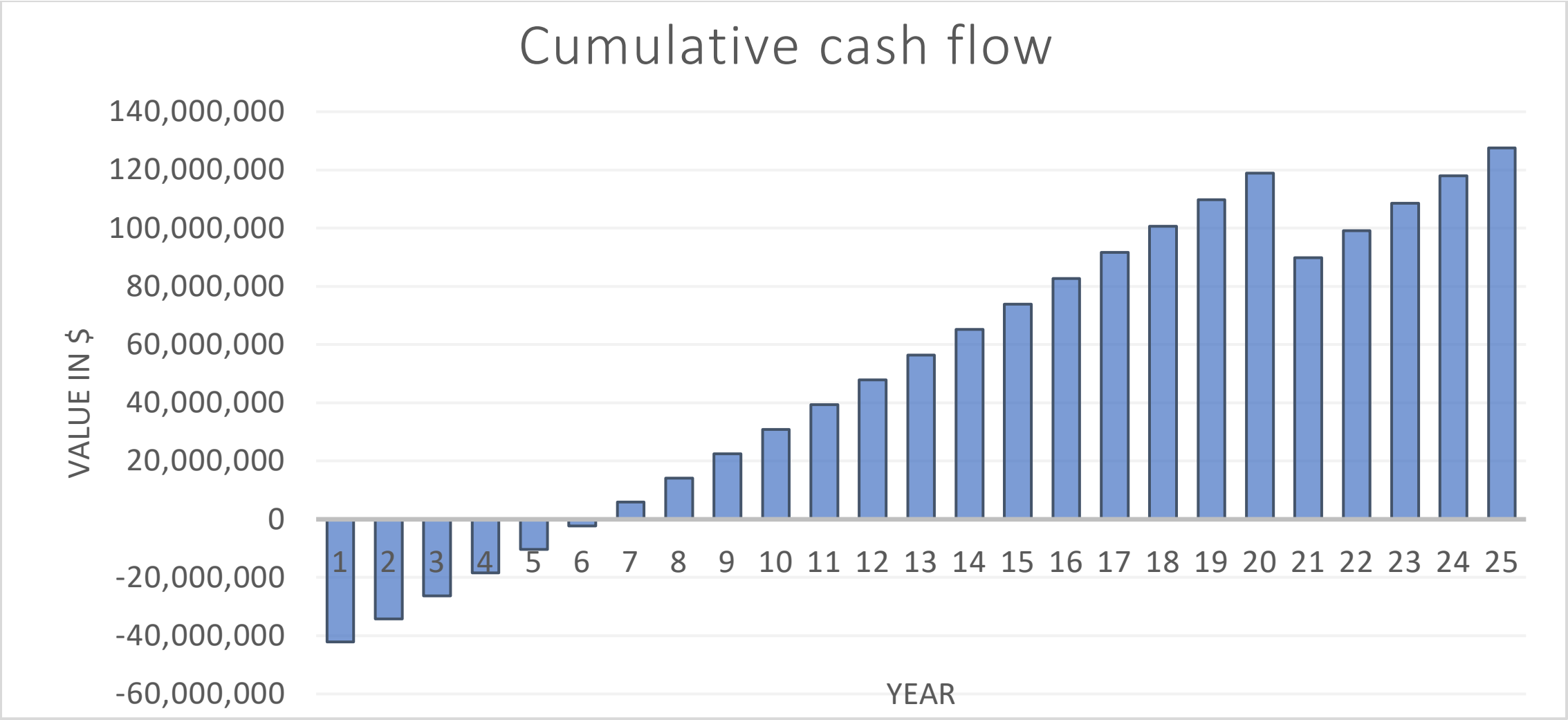

**Fig. 34: Cumulative cash flow of the project.**

The initial investment for implementing the project is evaluated at $42.1 million, and the total operating and maintenance cost of the whole system for 25 years is $7.46 million. It was found that the PV system and the set of converters each cost initially 4% of the initial investment while the rest of the cost is attributed to the wind turbines.

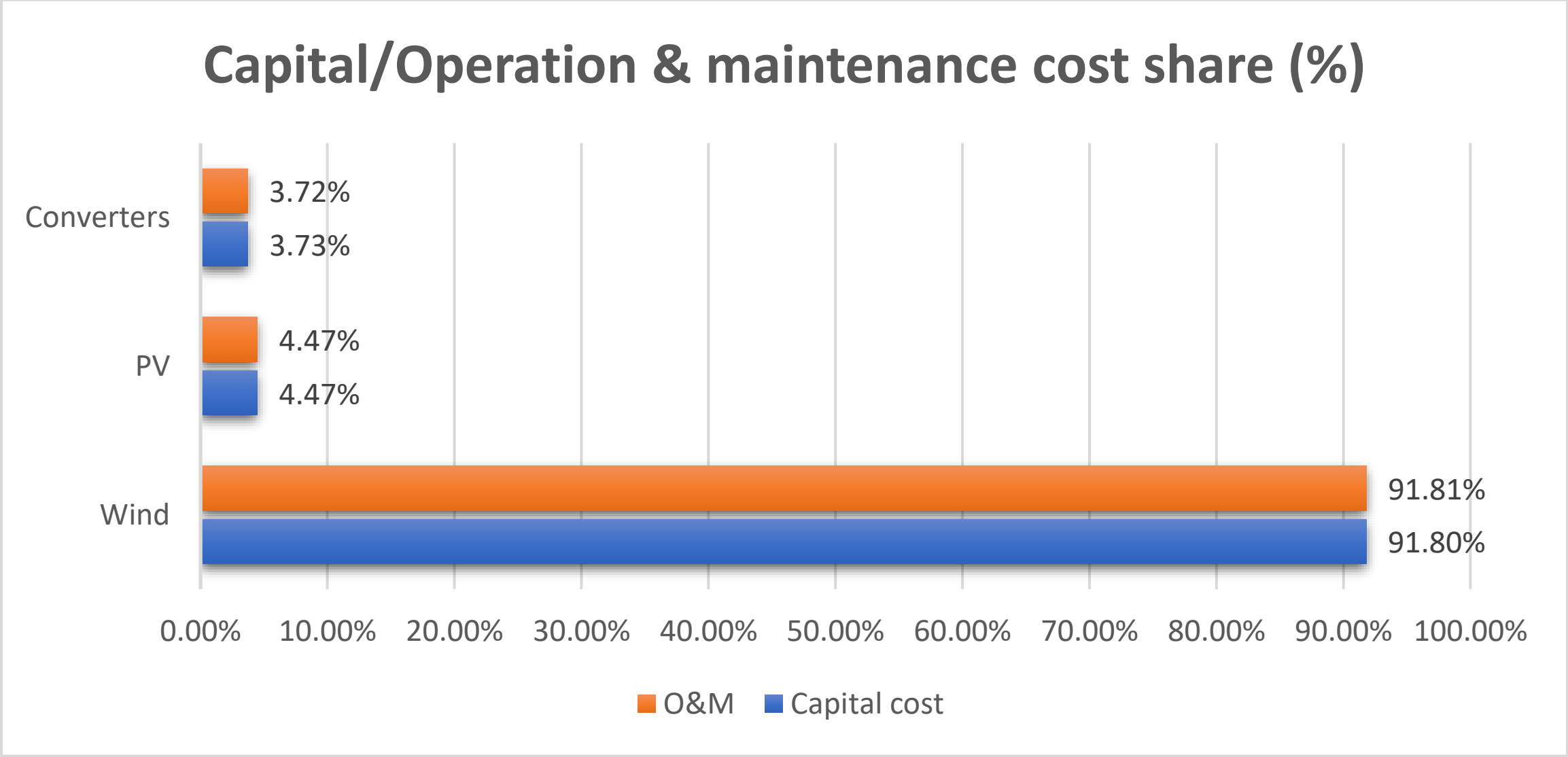

**Fig. 35: Components share of the total cost.**

The wind system dominates the share of both capital cost and operating and maintenance costs as seen in Fig. 35. Both PV system and converts account for 4% each of the total capital cost compared to 92% for the wind system. For every year of the project lifetime, the operating and maintenance of the PV array and converters require yearly around $15036 and $12530 respectively. On the other hand, the wind turbines require $309000 representing 91.8% of the total cost. Furthermore, the PV system produces 14.5% of to the total power generation versus 85.5% coming from the wind system. The monthly power production of both systems is shown in Fig. 36. There is a big difference in energy contribution that depends on two main factors; the different size of each system and the difference in the respective energy potential over the year.

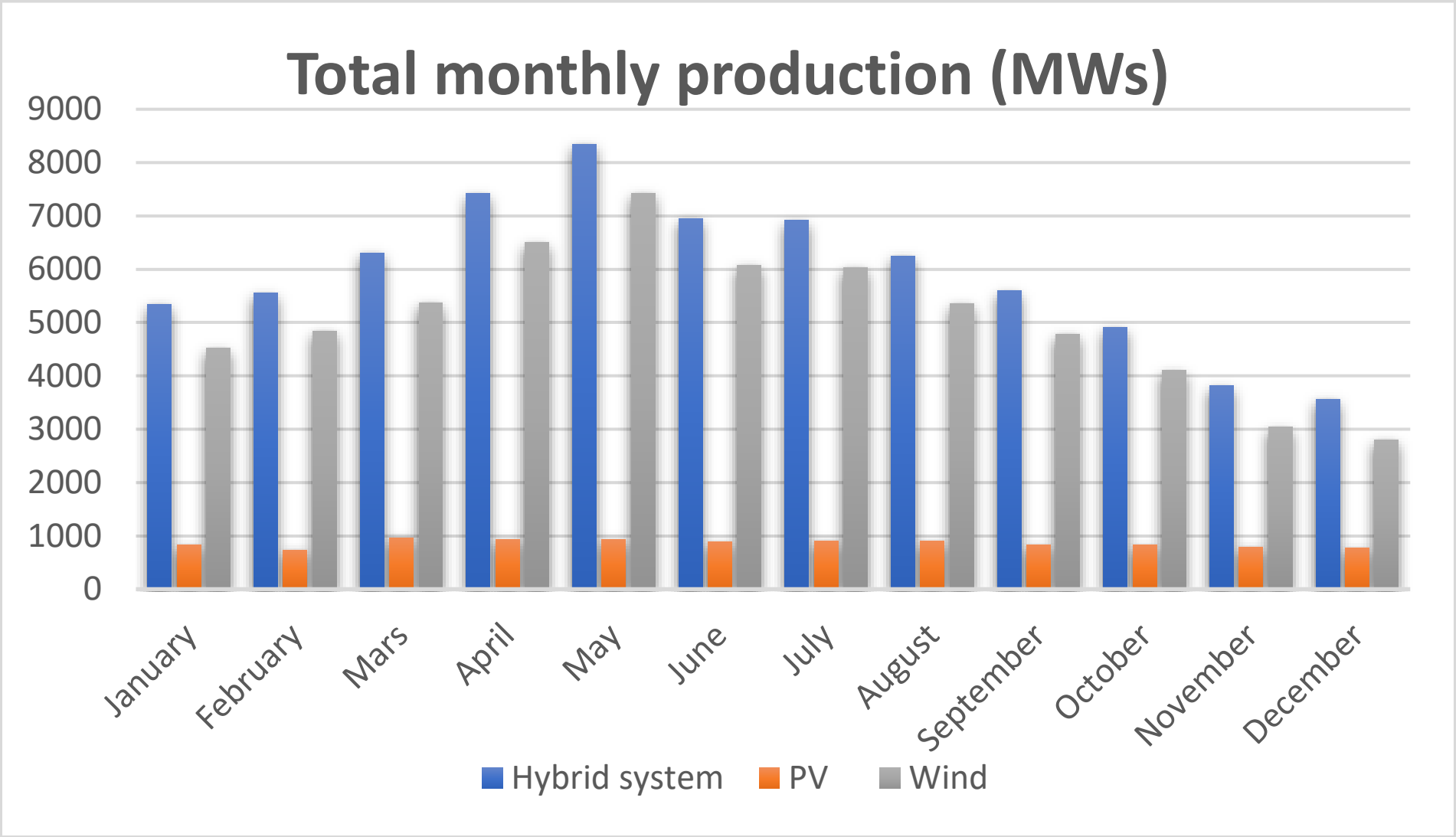

**Fig. 36: Monthly power production.**

## 3.9. Environmental Impact

From the environmental aspect, the designed system has a significant positive impact on the environment over the conventional power plant. It was found that implementing the hybrid system yields a remarkable reduction in harmful GHG emissions including $CO_2$, $SO_2$ and $NO_x$. For instance, a general gas-based combined cycle power plant produces around 400 kg/MWh. Therefore, integrating the system to the grid would dramatically reduce the amount of emitted GHG given the fact that electricity production in Algeria is mainly based on natural gas. Such emission reductions contribute to mitigating global warming and climate change. The summary of the amounts of avoided emissions have been tabulated below.

**Table 5: Reduced GHG emissions**

| Gas | Emissions reduction (kg/year) |
|---|---|
| Carbon Dioxide ($CO_2$) | 28 393 342 |
| Sulfur Dioxide ($SO_2$) | 194 494 |
| Nitrogen Oxides ($NO_x$) | 95 118 |

# 4. Conclusion

This article has addressed the design, modeling and control of a large-scale hybrid PV-wind grid connected system. The developed system has been tested for the Adrar region situated in the Algerian desert due to its relevant wind and solar energy resources. An innovative control approach using optimized PSO-PI controllers has been proposed to control the different inverters and converters of the wind turbines and PV system, where a new improved PSO algorithm is implemented to enhance the performance of the controllers. The performance of the designed system was assessed using hourly weather data for January and August. These two months represent the critical periods of winter and summer where the energy demand and energy resources vary the most. To evaluated the contribution of the system to the grid, its generated power was compared to the one generated by the central power plant of the region in question. It was found that the system is able to cover an important share of the total power consumption in the region but it mostly depends on energy resource variations. Furthermore, an economic feasibility and environmental impact of the hybrid system was carried out using HOMER. The results showed that the developed system proved to be technically and economically feasible. From an economical point of view, the project is cost-effective, with short a pay-back time and significant long-term profits. The environmental analysis revealed that implementing this project yields a significant yearly reduction in GHG emissions. The proposed hybrid system is an adequate solution to power shortages and grid problems faced in the region of Adrar during hot seasons. The proposed solution falls in line with the plan of Algeria to integrate wind and solar energy in its energy mix by 2030.

## DECLARATION

The authors emphasize that there are no conflicts of interest about the publication of this article.

## ACKNOWLEDGEMENTS

The authors would like to express their gratitude to the African Union and the Pan African University Institute of Water and Energy Sciences (including climate change) for their financial support to conduct this research

# Appendix

## A1

**Table 2: Main characteristics of the AS-6M module**

| | | |
|---|---|---|
| $I_{sc}$ | Short circuit current (A) | 9.68 |
| $k_i$ | Short circuit current of the cell at 25 °C and 1000 w/m$^2$ | $5 \times 10^{-4}$ |
| q | Electron charge (C) | $1.6 \times 10^{-19}$ |
| $V_{oc}$ | Open circuit voltage (V) | 46.5 |
| n | The ideality factor of the diode | 1.2 |
| K | Boltzmann’s constant (J/K) | $1.38 \times 10^{-23}$ |
| $E_{g0}$ | Bandgap energy of the semiconductor (eV) | 1.12 |
| $N_s$ | Number of cells connected in series | 72 |
| $R_s$ | Series resistor (Ω) | 0.100 |
| $R_p$ | Shunt resistor (Ω) | 515 |
| $P_{max}$ | Maximum output power(W) | 350 |
| $V_{mpp}$ | Maximum power voltage(V) | 37.9 |
| $I_{mpp}$ | Maximum power current(A) | 9.24 |

## A2

**Table 3: Power coefficients**

| $C_1$ | $C_2$ | $C_3$ | $C_4$ | $C_5$ | $C_6$ |
|---|---|---|---|---|---|
| 0.5176 | 116 | 0.4 | 5 | 21 | 0.0068 |

**A3**

**Table 4: Parameters of the wind system**

| **System** | **Parameter** |
|---|---|
| Turbine | Number of blades:3<br>R=45m<br>$J_t$=1.4 .$10^6$ kg.m²<br>$V_n$=13 m/s , $N_{tn}$=19 rpm |
| Gearbox | G=100 |
| DFIG | $U_r = U_s$= 690 V<br>$P_n$ = 3MW, f = 50 Hz, p = 2, Nr/Ns ≈ 1<br>$R_s$ = 2.97 mΩ, $R_r$ = 3.82 mΩ<br>$L_{fs}$ = 121 µH, $L_{fr}$ = 57.3 µH, $L_m$ = 12.12 mH<br>$J_g$ = 114 kg.m² |
| DC bus | C = 38 mF, $V_{dc}$ = 1200 V |
| RL Filter | $R_f$ = 0.75 Ω, $L_f$ = 0.75 mH |
| Grid | U = 690 V, f = 50 Hz |
| MPPT | $\lambda_{opt}$ = 7.07<br>$C_{pmax}$ = 0.35 |
| Rotor currents | Response time: $t_r$ = 0.05 s<br>$K_{pr} = 3.\sigma.\frac{L_r}{t_r} \approx 0.01062$<br>$K_{ir} = 3.\frac{R_r}{t_r} \approx 0.2292$ |
| Stator currents | Response time: $t_{rf}$ = 0.01 s<br>$K_{pf} = 3.\frac{L_f}{t_{rf}} \approx 0.3$<br>$K_{if} = 3.\frac{R_{rf}}{t_{rf}} \approx 30$ |